\documentclass[10pt,twoside]{article}
\usepackage{amsmath,amssymb,amsthm,latexsym} 
\usepackage{amsfonts,array}

\def\cal{\mathcal}

\def\dim{\mathop{\rm dim}\nolimits}

\def\tr{\mathop{\rm tr}\nolimits}
\def\Hom{\mathop{\rm Hom}\nolimits}

\def\id{\mathop{\rm id}\nolimits}

\def\K{F}

\def\gl{\mathfrak{gl}}
\def\so{\mathfrak{so}}

\def\V{{\cal V}}

\def\id{\rm id}

\def\gl{{\mathfrak{gl}}}
\def\so{{\mathfrak{so}}}
\def\card{{\rm card}}

\def\Hom{{\rm Hom}}
\def\Brk{{\bf Br}_k}
\def\Brt{{\bf Br}_2}
\def\Br3{{\bf Br}_3}
\def\bino(#1,#2){\binom{#1}{#2}}
\def\Du{u}
\def\Dx{x}

\def\id{{\rm id}}
\def\tr{{\rm tr}}
\def\ord{{\rm ord}}
\def\vspan{{\rm span}}
\def\mod{{\rm mod}}
\def\wc{{\tilde{c}}}

\def\Wgl{{\overline{W}_\gl}}
\def\Wso{{\overline{W}_\so}}
\def\V{{\cal V}} 
\def\Vnl{\V_{n,\ell}} 
\def\Ab{{\bar{\cal A}}} 
\def\Anl{{\cal A}_{n,\ell}} 
\def\Abnl{{\Ab_{n,\ell}}} 
\def\Abn1{{\Ab_n}} 
\def\P{{\cal P}} 
\def\Pn{\P_n} 
\def\Hvnl{{{\cal H}_{n,\ell}}} 
\def\Fvnl{{{\cal F}_{n,\ell}}} 
\def\Hvn1{{{\cal H}_{n}}} 
\def\Fvn1{{{\cal F}_{n}}} 
\def\Hvne{{{\cal H}_{n,1}}} 
\def\Fvne{{{\cal F}_{n,1}}} 
\def\Fvm1{{{\cal F}_{m,1}}} 
\def\Hnl{{{\cal H}_{n,\ell}'}} 
\def\Fnl{{{\cal F}_{n,\ell}'}} 
\def\Fil{{{\cal F}_{i,\ell}'}} 
\def\Fnmil{{{\cal F}_{n-i,\ell}'}}

\def\Fvi1{{{\cal F}_{i}}} 
\def\Fvnmi1{{{\cal F}_{n-i}}} 
\def\Knl{{{\cal F}_{n,\ell}'}} 
\def\Hn1{{{\cal H}_n'}} 
\def\Kn1{{{\cal F}_n'}} 
\def\Fn1{{{\cal F}_n'}} 
\def\H61{{{\cal H}_6'}} 
\def\K61{{{\cal F}_6'}}

\def\btrace{
\thicklines
\qbezier[100](-1,-1)(-1,-1.9)(1,-1.9)\qbezier[100](1,-1.9)(3,-1.9)(3,-1)
\qbezier[80](0,-1)(0,-1.6)(1,-1.6)\qbezier[80](1,-1.6)(2.5,-1.6)(2.5,-1)
\qbezier[60](1,-1)(1,-1.3)(1.5,-1.3)\qbezier[60](1.5,-1.3)(2,-1.3)(2,-1)
\put(2,-1){\line(0,1){2}}
\put(2.5,-1){\line(0,1){2}}
\put(3,-1){\line(0,1){2}}
\qbezier[100](-1,1)(-1,1.9)(1,1.9)\qbezier[100](1,1.9)(3,1.9)(3,1)
\qbezier[80](0,1)(0,1.6)(1,1.6)\qbezier[80](1,1.6)(2.5,1.6)(2.5,1)
\qbezier[60](1,1)(1,1.3)(1.5,1.3)\qbezier[60](1.5,1.3)(2,1.3)(2,1)
}

\def\bline[#1,#2]{
\thicklines
\ifnum#1>#2\bline[#2,#1]
\else\ifnum#1<#2
\ifcase#1
\ifcase#2
\or\qbezier[60](-1,-1)(-1,-0.5)(-0.5,-0.5)
\qbezier[60](-0.5,-0.5)(0,-0.5)(0,-1)
\or\qbezier[80](-1,-1),(-1,-0.5)(0,-0.5)
\qbezier[80](0,-0.5)(1,-0.5)(1,-1)
\or\put(-1,-1){\line(1,1){2}}
\or\put(-1,-1){\line(1,2){1}}
\or\put(-1,-1){\line(0,1){2}}
\fi
\or\ifcase#2\or
\or\qbezier[80](0,-1),(0,-0.5)(0.5,-0.5)
\qbezier[80](0.5,-0.5)(1,-0.5)(1,-1)
\or\put(0,-1){\line(1,2){1}}
\or\put(0,-1){\line(0,1){2}}
\or\put(0,-1){\line(-1,2){1}}
\fi
\or\ifcase#2\or\or
\or\put(1,-1){\line(0,1){2}}
\or\put(1,-1){\line(-1,2){1}}
\or\put(1,-1){\line(-1,1){2}}
\fi
\or\ifcase#2\or\or\or
\or\qbezier[60](1,1)(1,0.5)(0.5,0.5)
\qbezier[60](0.5,0.5)(0,0.5)(0,1)
\or\qbezier[80](1,1),(1,0.5)(0,0.5)
\qbezier[80](0,0.5)(-1,0.5)(-1,1)
\fi
\or\ifcase#2\or\or\or\or
\or\qbezier[60](-1,1)(-1,0.5)(-0.5,0.5)
\qbezier[60](-0.5,0.5)(0,0.5)(0,1)
\fi\fi\fi\fi
}

\theoremstyle{plain}
\newtheorem{lemma}{Lemma}%[section]
\newtheorem{prop}[lemma]{Proposition}
\newtheorem{theorem}[lemma]{Theorem}
\newtheorem{coro}[lemma]{Corollary}
\theoremstyle{definition}

\theoremstyle{remark}

\newcommand{\Z}{{\ensuremath{\mathbb Z}}}
\newcommand{\Q}{{\ensuremath{\mathbb Q}}}
\newcommand{\R}{{\ensuremath{\mathbb R}}}

\def\DottedCircle{
\bezier{4}(0.966,-0.259)(1.04,0)(0.966,0.259)
\bezier{4}(0.966,0.259)(0.897,0.518)(0.707,0.707)
\bezier{4}(0.707,0.707)(0.518,0.897)(0.259,0.966)
\bezier{4}(0.259,0.966)(0,1.04)(-0.259,0.966)
\bezier{4}(-0.259,0.966)(-0.518,0.897)(-0.707,0.707)
\bezier{4}(-0.707,0.707)(-0.897,0.518)(-0.966,0.259)
\bezier{4}(-0.966,0.259)(-1.04,0)(-0.966,-0.259)
\bezier{4}(-0.966,-0.259)(-0.897,-0.518)(-0.707,-0.707)
\bezier{4}(-0.707,-0.707)(-0.518,-0.897)(-0.259,-0.966)
\bezier{4}(-0.259,-0.966)(0,-1.04)(0.259,-0.966)
\bezier{4}(0.259,-0.966)(0.518,-0.897)(0.707,-0.707)
\bezier{4}(0.707,-0.707)(0.897,-0.518)(0.966,-0.259)
}
\def\FullCircle{
\thicklines
\put(0,0){\circle{2}}
}
\def\Endpoint[#1]{
}
\def\Arc[#1]{
\thicklines                     % this can be changed!
\ifcase#1
\bezier{25}(0.966,-0.259)(1.04,0)(0.966,0.259)
\or
\bezier{25}(0.966,0.259)(0.897,0.518)(0.707,0.707)
\or
\bezier{25}(0.707,0.707)(0.518,0.897)(0.259,0.966)
\or
\bezier{25}(0.259,0.966)(0,1.04)(-0.259,0.966)
\or
\bezier{25}(-0.259,0.966)(-0.518,0.897)(-0.707,0.707)
\or
\bezier{25}(-0.707,0.707)(-0.897,0.518)(-0.966,0.259)
\or
\bezier{25}(-0.966,0.259)(-1.04,0)(-0.966,-0.259)
\or
\bezier{25}(-0.966,-0.259)(-0.897,-0.518)(-0.707,-0.707)
\or
\bezier{25}(-0.707,-0.707)(-0.518,-0.897)(-0.259,-0.966)
\or
\bezier{25}(-0.259,-0.966)(0,-1.04)(0.259,-0.966)
\or
\bezier{25}(0.259,-0.966)(0.518,-0.897)(0.707,-0.707)
\or
\bezier{25}(0.707,-0.707)(0.897,-0.518)(0.966,-0.259)
\fi}
\def\DottedArc[#1]{
\ifcase#1
\bezier{4}(0.966,-0.259)(1.04,0)(0.966,0.259)
\or
\bezier{4}(0.966,0.259)(0.897,0.518)(0.707,0.707)
\or
\bezier{4}(0.707,0.707)(0.518,0.897)(0.259,0.966)
\or
\bezier{4}(0.259,0.966)(0,1.04)(-0.259,0.966)
\or
\bezier{4}(-0.259,0.966)(-0.518,0.897)(-0.707,0.707)
\or
\bezier{4}(-0.707,0.707)(-0.897,0.518)(-0.966,0.259)
\or
\bezier{4}(-0.966,0.259)(-1.04,0)(-0.966,-0.259)
\or
\bezier{4}(-0.966,-0.259)(-0.897,-0.518)(-0.707,-0.707)
\or
\bezier{4}(-0.707,-0.707)(-0.518,-0.897)(-0.259,-0.966)
\or
\bezier{4}(-0.259,-0.966)(0,-1.04)(0.259,-0.966)
\or
\bezier{4}(0.259,-0.966)(0.518,-0.897)(0.707,-0.707)
\or
\bezier{4}(0.707,-0.707)(0.897,-0.518)(0.966,-0.259)
\fi}
\def\Chord[#1,#2]{
\thinlines
\ifnum#1>#2\Chord[#2,#1]
\else\ifnum#1<#2
\ifcase#1
\ifcase#2
\or\qbezier(1,0)(0.516,0.138)(0.866,0.5)
\or\qbezier(1,0)(0.45,0.26)(0.5,0.866)
\or\qbezier(1,0)(0.327,0.327)(0,1)
\or\qbezier(1,0)(0.179,0.311)(-0.5,0.866)
\or\qbezier(1,0)(0.0536,0.2)(-0.866,0.5)
\or\put(1, 0){\line(-2, 0){2}}
\or\qbezier(1,0)(0.0536,-0.2)(-0.866,-0.5)
\or\qbezier(1,0)(0.179,-0.311)(-0.5,-0.866)
\or\qbezier(1,0)(0.327,-0.327)(0,-1)
\or\qbezier(1,0)(0.45,-0.26)(0.5,-0.866)
\or\qbezier(1,0)(0.516,-0.138)(0.866,-0.5)
\fi
\or\ifcase#2\or
\or\qbezier(0.866,0.5)(0.378,0.378)(0.5,0.866)
\or\qbezier(0.866,0.5)(0.26,0.45)(0,1)
\or\qbezier(0.866,0.5)(0.12,0.446)(-0.5,0.866)
\or\qbezier(0.866,0.5)(0,0.359)(-0.866,0.5)
\or\qbezier(0.866,0.5)(-0.0536,0.2)(-1,0)
\or\put(0.866, 0.5){\line(-5, -3){1.73}}
\or\qbezier(0.866,0.5)(0.146,-0.146)(-0.5,-0.866)
\or\qbezier(0.866,0.5)(0.311,-0.179)(0,-1)
\or\qbezier(0.866,0.5)(0.446,-0.12)(0.5,-0.866)
\or\qbezier(0.866,0.5)(0.52,0)(0.866,-0.5)
\fi
\or\ifcase#2\or\or
\or\qbezier(0.5,0.866)(0.138,0.516)(0,1)
\or\qbezier(0.5,0.866)(0,0.52)(-0.5,0.866)
\or\qbezier(0.5,0.866)(-0.12,0.446)(-0.866,0.5)
\or\qbezier(0.5,0.866)(-0.179,0.311)(-1,0)
\or\qbezier(0.5,0.866)(-0.146,0.146)(-0.866,-0.5)
\or\put(0.5, 0.866){\line(-3, -5){1}}
\or\qbezier(0.5,0.866)(0.2,-0.0536)(0,-1)
\or\qbezier(0.5,0.866)(0.359,0)(0.5,-0.866)
\or\qbezier(0.5,0.866)(0.446,0.12)(0.866,-0.5)
\fi
\or\ifcase#2\or\or\or
\or\qbezier(0,1.)(-0.138,0.516)(-0.5,0.866)
\or\qbezier(0,1.)(-0.26,0.45)(-0.866,0.5)
\or\qbezier(0,1.)(-0.327,0.327)(-1,0)
\or\qbezier(0,1.)(-0.311,0.179)(-0.866,-0.5)
\or\qbezier(0,1.)(-0.2,0.0536)(-0.5,-0.866)
\or\put(0, 1){\line(0, -2){2}}
\or\qbezier(0,1.)(0.2,0.0536)(0.5,-0.866)
\or\qbezier(0,1.)(0.311,0.179)(0.866,-0.5)
\fi
\or\ifcase#2\or\or\or\or
\or\qbezier(-0.5,0.866)(-0.378,0.378)(-0.866,0.5)
\or\qbezier(-0.5,0.866)(-0.45,0.26)(-1,0)
\or\qbezier(-0.5,0.866)(-0.446,0.12)(-0.866,-0.5)
\or\qbezier(-0.5,0.866)(-0.359,0)(-0.5,-0.866)
\or\qbezier(-0.5,0.866)(-0.2,-0.0536)(0,-1)
\or\put(-0.5, 0.866){\line(3, -5){1}}
\or\qbezier(-0.5,0.866)(0.146,0.146)(0.866,-0.5)
\fi
\or\ifcase#2\or\or\or\or\or
\or\qbezier(-0.866,0.5)(-0.516,0.138)(-1,0)
\or\qbezier(-0.866,0.5)(-0.52,0)(-0.866,-0.5)
\or\qbezier(-0.866,0.5)(-0.446,-0.12)(-0.5,-0.866)
\or\qbezier(-0.866,0.5)(-0.311,-0.179)(0,-1)
\or\qbezier(-0.866,0.5)(-0.146,-0.146)(0.5,-0.866)
\or\put(-0.866, 0.5){\line(5, -3){1.73}}
\fi
\or\ifcase#2\or\or\or\or\or\or
\or\qbezier(-1,0)(-0.516,-0.138)(-0.866,-0.5)
\or\qbezier(-1,0)(-0.45,-0.26)(-0.5,-0.866)
\or\qbezier(-1,0)(-0.327,-0.327)(0,-1)
\or\qbezier(-1,0)(-0.179,-0.311)(0.5,-0.866)
\or\qbezier(-1,0)(-0.0536,-0.2)(0.866,-0.5)
\fi
\or\ifcase#2\or\or\or\or\or\or\or
\or\qbezier(-0.866,-0.5)(-0.378,-0.378)(-0.5,-0.866)
\or\qbezier(-0.866,-0.5)(-0.26,-0.45)(0,-1)
\or\qbezier(-0.866,-0.5)(-0.12,-0.446)(0.5,-0.866)
\or\qbezier(-0.866,-0.5)(0,-0.359)(0.866,-0.5)
\fi
\or\ifcase#2\or\or\or\or\or\or\or\or
\or\qbezier(-0.5,-0.866)(-0.138,-0.516)(0,-1)
\or\qbezier(-0.5,-0.866)(0,-0.52)(0.5,-0.866)
\or\qbezier(-0.5,-0.866)(0.12,-0.446)(0.866,-0.5)
\fi
\or\ifcase#2\or\or\or\or\or\or\or\or\or
\or\qbezier(0,-1.)(0.138,-0.516)(0.5,-0.866)
\or\qbezier(0,-1.)(0.26,-0.45)(0.866,-0.5)
\fi
\or\ifcase#2\or\or\or\or\or\or\or\or\or\or
\or\qbezier(0.5,-0.866)(0.378,-0.378)(0.866,-0.5)
\fi\fi\fi\fi}
\def\FullChord[#1,#2]{
\Endpoint[#1]
\Endpoint[#2]
\Arc[#1]
\Arc[#2]
\Chord[#1,#2]
}
\def\EndChord[#1,#2]{
\Endpoint[#1]
\Endpoint[#2]
\Chord[#1,#2]
}
\def\Picture#1{
\begin{picture}(2,1)(-1,-0.167) 
#1
\end{picture}
}
\def\DottedChordDiagram[#1,#2]{
\Picture{\DottedCircle \FullChord[#1,#2]}
}
\def\ExtChord[#1,#2]{
\Endpoint[#1]\Endpoint[#2]
\thinlines
\ifnum#1>#2\ExtChord[#2,#1]
\else\ifnum#1<#2
\ifcase#1
\or\ifcase#2
\or\or\or\or\or\or\or\or\or
\qbezier[80](0,-1)(-0.1,-1.4)(0.25,-1.35)
\qbezier[80](0.25,-1.35)(1.35,-1.2)(1.35,0)
\qbezier[80](1.35,0)(1.35,0.95)(0.866,0.5)
\fi
\or\or\ifcase#2
\or\or\or\or\or\or\or\or\or\or\or
\qbezier[80](0,1)(-0.1,1.4)(0.25,1.35)
\qbezier[80](0.25,1.35)(1.35,1.2)(1.35,0)
\qbezier[80](1.35,0)(1.35,-0.95)(0.866,-0.5)
\fi
\or\ifcase#2
\or\or\or\or\or\or\or\or
\qbezier[80](-0.5,0.866)(-0.65,1.1)(-0.75,1.1)
\qbezier[60](-0.75,1.1)(-1.35,1.1)(-1.35,0)
\qbezier[60](-0.5,-0.866)(-0.65,-1.1)(-0.75,-1.1)
\qbezier[80](-0.75,-1.1)(-1.35,-1.1)(-1.35,0)
\fi\fi\fi\fi
}

\def\YEAR{\year}\newcount\VOL\VOL=\YEAR\advance\VOL by-1995
\def\firstpage{275}\def\lastpage{298}
\def\received{February 17, 2000}\def\revised{}
\def\communicated{G\"unter M. Ziegler}

\makeatletter
\def\magnification{\afterassignment\m@g\count@}
\def\m@g{\mag=\count@\hsize6.5truein\vsize8.9truein\dimen\footins8truein}
\makeatother

\font\eightrm=cmr8
\font\caps=cmcsc10                    % Theorem, Lemma etc
\font\Caps=cccsc10 scaled \magstep1   % Title
\font\scaps=cmcsc8

\makeatletter
\@specialpagefalse\headheight=8.5pt
\def\DocMath{{\def\th{\thinspace}\scaps Documenta Math.}}
\renewcommand{\@oddfoot}{\hfill\scaps Documenta Mathematica 
    \number\VOL\  (\number\YEAR) \number\firstpage--\lastpage\hfill}
\renewcommand{\@evenfoot}{\ifnum\thepage>\lastpage\hfill\scaps
    Documenta Mathematica \number\VOL\  (\number\YEAR)\hfill\else\@oddfoot\fi}%

\renewcommand{\@evenhead}{%
    \ifnum\thepage>\lastpage\rlap{\thepage}\hfill%
    \else\rlap{\thepage}\slshape\leftmark\hfill\caps\SAuthor\hfill\fi}%
\renewcommand{\@oddhead}{%
    \ifnum\thepage=\firstpage{\DocMath\hfill\llap{\thepage}}%
    \else{\slshape\rightmark}\hfill\caps\STitle\hfill\llap{\thepage}\fi}%
\makeatother

\def\TSkip{\bigskip}
\newbox\TheTitle{\obeylines\gdef\GetTitle #1
\ShortTitle  #2
\SubTitle    #3
\Author      #4
\ShortAuthor #5
\EndTitle
{\setbox\TheTitle=\vbox{\baselineskip=20pt\let\par=\cr\obeylines%
\halign{\centerline{\Caps##}\cr\noalign{\medskip}\cr#1\cr}}%
	\copy\TheTitle\TSkip\TSkip%
\def\next{#2}\ifx\next\empty\gdef\STitle{#1}\else\gdef\STitle{#2}\fi%
\def\next{#3}\ifx\next\empty%
    \else\setbox\TheTitle=\vbox{\baselineskip=20pt\let\par=\cr\obeylines%
    \halign{\centerline{\caps##} #3\cr}}\copy\TheTitle\TSkip\TSkip\fi%
\centerline{\caps #4}\TSkip\TSkip%
\def\next{#5}\ifx\next\empty\gdef\SAuthor{#4}\else\gdef\SAuthor{#5}\fi%
\ifx\received\empty\relax
    \else\centerline{\eightrm Received: \received}\fi%
\ifx\revised\empty\TSkip%
    \else\centerline{\eightrm Revised: \revised}\TSkip\fi%
\ifx\communicated\empty\relax
    \else\centerline{\eightrm Communicated by \communicated}\fi\TSkip\TSkip%
\catcode'015=5}}\def\Title{\obeylines\GetTitle}
\def\Abstract{\begingroup\narrower
    \parskip=\medskipamount\parindent=0pt{\caps Abstract. }}
\def\EndAbstract{\par\endgroup\TSkip}

\long\def\MSC#1\EndMSC{\def\arg{#1}\ifx\arg\empty\relax\else
     {\par\narrower\noindent%
     1991 Mathematics Subject Classification: #1\par}\fi}

\long\def\KEY#1\EndKEY{\def\arg{#1}\ifx\arg\empty\relax\else
	{\par\narrower\noindent Keywords and Phrases: #1\par}\fi\TSkip}

\newbox\TheAdd\def\Addresses{\vfill\copy\TheAdd\vfill
    \ifodd\number\lastpage\vfill\eject\phantom{.}\vfill\eject\fi}
{\obeylines\gdef\GetAddress #1
\Address #2 
\Address #3
\Address #4
\EndAddress
{\def\xs{4.3truecm}\parindent=0pt
\setbox0=\vtop{{\obeylines\hsize=\xs#1\par}}\def\next{#2}
\ifx\next\empty % 1 address
     \setbox\TheAdd=\hbox to\hsize{\hfill\copy0\hfill}
\else\setbox1=\vtop{{\obeylines\hsize=\xs#2\par}}\def\next{#3}
\ifx\next\empty % 2 addresses
     \setbox\TheAdd=\hbox to\hsize{\hfill\copy0\hfill\copy1\hfill}
\else\setbox2=\vtop{{\obeylines\hsize=\xs#3\par}}\def\next{#4}
\ifx\next\empty\ % 3 addresses
     \setbox\TheAdd=\vtop{\hbox to\hsize{\hfill\copy0\hfill\copy1\hfill}
                \vskip20pt\hbox to\hsize{\hfill\copy2\hfill}}
\else\setbox3=\vtop{{\obeylines\hsize=\xs#4\par}}
     \setbox\TheAdd=\vtop{\hbox to\hsize{\hfill\copy0\hfill\copy1\hfill}
	        \vskip20pt\hbox to\hsize{\hfill\copy2\hfill\copy3\hfill}}
\fi\fi\fi\catcode'015=5}}\gdef\Address{\obeylines\GetAddress}

\begin{document} 
%%%% ------------- fill in your data below this line  -------------------

\Title 
The Number of 
Independent Vassiliev Invariants in 
the Homfly and Kauffman Polynomials
\ShortTitle 
Number of Independent Vassiliev Invariants in $H$ and $F$
\SubTitle   
\Author 
Jens Lieberum
\ShortAuthor
\EndTitle

\Abstract We consider vector spaces~$\Hvnl$ and~$\Fvnl$ spanned by 
the degree-$n$ coefficients in power series forms of the Homfly 
and Kauffman polynomials of links with~$\ell$ components. 
Generalizing previously known formulas, we determine the 
dimensions of the spaces~$\Hvnl$, $\Fvnl$ and $\Hvnl+\Fvnl$ for 
all values of~$n$ and~$\ell$. Furthermore, we show that for knots 
the algebra generated by $\bigoplus_n \Hvne+\Fvne$ is a polynomial 
algebra with~$\dim(\Hvne+\Fvne)-1=n+[n/2]-4$ generators in 
degree~$n\geq 4$ and one generator in degrees~$2$ and~$3$. 

\EndAbstract 
\MSC 
%%%%%    2000 Mathematics Subject Classification: 
57M25.
\EndMSC
\KEY 
%%%%%    Keywords and Phrases: 
Vassiliev invariants, link polynomials, 
Brauer algebra, Vogel´s algebra, dimensions.
\EndKEY
\Address Jens Lieberum 
MSRI 
1000 Centennial Drive 
Berkeley,~CA~94720-5070 
USA 
lieberum@msri.org 
\Address 
\Address 
\Address 
\EndAddress 
%
% Make sure the last tex command in your manuscript
% before the first \end or \bye is the command  \Addresses
%
%---------------------Here the prologue ends---------------------------------
%--------------------Here the manuscript starts------------------------------

%\input Hkn.tex
{\parindent0cm
\def\Upsilon{{Y}}
%\markright{The common specializations of the Homfly and Kauffman
%polynomials}
%\chapter{The common specializations of the Homfly and Kauffman
%polynomials}
%\input{Hkt.tex}

%We calculate the dimensions of the space of Vassiliev invariants
%coming from the Homfly polynomial of links, and of the space of Vassiliev
%invariants coming from the Kauffman polynomial of links.
%We show that the intersection of these spaces
%is spanned by the Vassiliev invariants coming from the
%Jones polynomial and from a polynomial called $\Upsilon$.
%We also show that linearly independent Vassiliev invariants of knots
%of degree $\geq 8$ coming from the Homfly or Kauffman polynomial are
%algebraically independent.

\section{Introduction}%\addcontentsline{toc}{section}{\numberline{}Introduction}

Soon after the discovery of the Jones polynomial~$V$ (\cite{Jon}),
two $2$-parameter
generali\-zations of it were introduced: the Homfly polynomial~$H$
(\cite{HOM}) and the Kauffman polynomial~$F$ (\cite{Ka2}) of oriented links.
Let $\Vnl$ be the vector space of $\Q$-valued Vassiliev invariants
of degree~$n$ of links with~$\ell$ components.
After a substitution of parameters,
the polynomial~$H$ (resp.~$F$)
can be written as a power series in an indeterminate~$h$, such that
the coefficient of~$h^n$ is a polynomial-valued Vassiliev invariant~$p_n$
(resp.~$q_n$) of degree~$n$.
Let~$\Hvnl$ (resp.~$\Fvnl$) be the vector space generated by the coefficients
of~$p_n$ (resp.~$q_n$) regarded as a subspace of~$\Vnl$.
The dimensions of~$\Hvnl$ and~$\Fvnl$ have been determined in
\cite{Men} for~$n\geq 0$ and~$\ell=1$ and partial results were also known for~$\ell>1$.
We complete these formulas by calculating $\dim\Hvnl$, $\dim\Fvnl$ and~$\dim(\Hvnl+\Fvnl)$
for~$n\geq 0$ and all pairs~$(n, \ell)$. 

\begin{theorem}\label{t:dimtab}
(1) For all $n,\ell\geq 1$ we have\nopagebreak{}

%\begin{equation}\label{dimhnl}
$$
\dim \Hvnl=\min\left\{n,\left[\frac{n-1+\ell}{2}\right]\right\}=
\left\{
\begin{array}{ll}
n & \mbox{if $n<\ell$,}\\
\left[\frac{n-1+\ell}{2}\right] & \mbox{if $n\geq \ell$.}
\end{array}
\right.
%\end{equation}
$$

(2) If $n\geq 4$, then\nopagebreak{}

%\begin{equation}\label{dimknl}
$$
\dim \Fvnl=\left\{
\begin{array}{ll}
n-1 & \mbox{if $\ell=1$,}\\
2n-1 & \mbox{if $\ell\geq 2$ and $n\leq \ell$,}\\
n+\ell-1 & \mbox{if $\ell\geq 2$ and $n\geq \ell$}.
\end{array}
\right.
%\end{equation}
$$

The values of $\dim\Fvnl$ for $n\leq 3$ are given in the following table\nopagebreak{}

$$
\begin{array}{|l|c|c|c|c|c|c|c|c|}
\hline
(n,\ell)                & (1,1) & (1,\geq 2) & (2,1) & (2,2) & (2,\geq 3) & (3,1) & (3,2) & (3,\geq 3)\\
\hline
\dim \Fvnl            & 0     & 1          & 1     & 2     & 3          & 1     & 4     & 5\\ 
\hline
\end{array}
$$

(3) For all $n,\ell\geq 1$ we have\nopagebreak{}

%\begin{equation}
$$
\dim(\Hvnl\cap \Fvnl)=\min\{\dim \Hvnl,2\}.%=\dim\left(\vspan\{r_n^\ell,
%y_n^\ell\}\right).
%\end{equation}
$$
\end{theorem}

%It was shown by H.\ Lamaugarny that, in a sense made precise in~\cite{Lam},
%all common specializations of the polynomials~$H$ and~$F$ can be expressed
%in terms of the Jones polynomial~$V$, the number of components of a link, and
%the linking number.
In the framework of Vassiliev invariants it is natural to consider
the elements of~$\bigoplus_{n,\ell} (\Hvnl\cap\Fvnl)$ as
{\em the common specializations of $H$ and $F$}.
It is known that a one-variable
polynomial~$\Upsilon$ (\cite{CoG}, \cite{Kn1}, \cite{Lik}, \cite{Lie}, \cite{Sul}) 
appears as a lowest coefficient in~$H$ and~$F$.
This is used in the proof of Theorem~\ref{t:dimtab} to derive lower bounds 
for~$\dim(\Hvnl\cap\Fvnl)$. Let~$r_n^\ell$ be the coefficient 
of~$h^n$ in the Jones polynomial~$V(e^{h/2})$ 
and let~$y_n^\ell$ be the coefficient of~$h^n$ in~$\Upsilon(e^{h/2})$.
Then
we have $r_n^\ell,y_n^\ell\in\Hvnl\cap\Fvnl$.
The following corollary to the proof of Theorem~\ref{t:dimtab} says that the Jones 
polynomial~$V$ and the polynomial~$Y$ are {\em the only} common specializations 
of~$H$ and~$F$ in the sense above (compare~\cite{Lam} for common specializations in a 
different sense).

\begin{coro}\label{c:hkjy}
For all $n\geq 0,\ell\geq 1$ we have
$\Hvnl\cap\Fvnl=\vspan\{r_n^\ell, y_n^\ell\}$.
\end{coro}

The main part of the proofs of Theorem~\ref{t:dimtab} and Corollary~\ref{c:hkjy}
will not be given on the level of link invariants,
but on the level of weight systems. A weight system of degree~$n$
is a linear form
on a space $\Abnl$ generated by certain trivalent graphs with~$\ell$ distinguished
oriented circles and~$2n$ vertices called {\em trivalent diagrams}.
There exists a surjective map~$W$ from~$\Vnl$ to the space
$\Abnl^*=\Hom(\Abnl,\Q)$ of weight systems.
The restriction of~$W$ to~$\Hvnl+\Fvnl$ is injective.
So we may study the spaces~$\Hnl=W(\Hvnl)$ and
$\Fnl=W(\Fvnl)\subseteq\Abnl^*$
instead of~$\Hvnl$ and~$\Fvnl$.
Using an explicit description of weight systems in~$\Hnl$ and~$\Fnl$ we
%It will be easy to 
derive upper bounds for
$\dim \Hnl$ and $\dim \Fnl$.
We obtain an upper bound for~$\dim(\Hnl+\Fnl)$ from a lower bound 
for~$\dim(\Hnl\cap\Fnl)$.
%Our main work will be to 
We evaluate 
the weight systems in~$\Hnl$ and~$\Fnl$
on many trivalent diagrams which gives us
%In this way we 
lower bounds for
$\dim \Hnl$, $\dim \Fnl$ and $\dim (\Hnl+\Fnl)$.
%Fortunately,
These lower bounds always coincide with the upper bounds.
The resulting dimension formulas %are of independent interest and
will imply Theorem~\ref{t:dimtab}.
%The use of weight systems in the proof has several advantages.

For simplicity of notation we will drop
the index~$\ell$ when~$\ell=1$.
The fact that the Jones polynomial and the square of the Jones polynomial
appear by choosing special values of parameters of the Kauffman polynomial
gives us quadratic relations between elements of~$\bigoplus_{n=0}^\infty \Fvnl$. 
We will use the Hopf algebra structure
of~$\Ab=\bigoplus_{n=0}^\infty \Ab_{n}$ to show that
%in degree~$n\geq 9$, linearly
we know
all algebraic relations between elements of $\bigoplus_{n=0}^\infty \Hvn1+\Fvn1$:
%The algebra generated by $\bigoplus_{n=0}^\infty \Hne+\Fne$ turns out to be
%a polynomial algebra with $\dim(\Hne+\Fne)-1=n+[n/2]-4$ generators in
%degree~$n\geq 4$ and one generator in degrees~$2$ and~$3$.

\begin{theorem}\label{t:alggen}
The algebra generated by
$
\bigoplus_{n=0}^\infty \Hvn1+\Fvn1
$
is a polynomial algebra with
 
$$\max\{\dim(\Hvn1+\Fvn1)-1,1\}=\max\{n+[n/2]-4,1\}$$ 

generators in degree~$n\geq 2$. 
\end{theorem}

If knot invariants $v_i$ satisfy $v_i(K_1)=v_i(K_2)$, then
polynomials in the invariants~$v_i$ also cannot
distinguish the knots~$K_1$ and~$K_2$. 
By Theorem~\ref{t:alggen} there is only one algebraic relation between elements
$v_i\in\bigoplus_{n=1}^{m-1} ({\cal H}_n + {\cal F}_n)$ and elements 
of~${\cal H}_m + {\cal F}_m$ in each degree $m \geq 4$. 
This gives us a hint 
why it is possible to distinguish many knots by
comparing their Homfly and Kauffman polynomials.

The plan of the paper is the following.
In Section~\ref{linkpolys} we recall the definitions of the
link polynomials $H$, $F$, $V$, $\Upsilon$, and we
give the exact definitions of~$\Hvnl$ and~$\Fvnl$.
%We also recall that the coefficients of $h^n$ in
%$\Upsilon(e^{h/2})$ and in $V(e^{h/2})$ are in $\Hvnl\cap\Fvnl$.
Then we express relations between these polynomials in terms of
Vassiliev invariants.
%elements of $\bigoplus_{n=0}^\infty(\Hvnl+\Fvnl)$.
In Section~\ref{dimtabs} we define~$\Abnl$
and recall the connection between the Vassiliev invariants in~$\Hvnl+\Fvnl$ and
their weight systems in~$\Hnl+\Fnl$. 
In Section~\ref{upperb} we use a direct combinatorial
description of the weight
systems in~$\Hnl$ and~$\Fnl$ to derive upper bounds for~$\dim\Hnl$ and~$\dim\Fnl$.
For the proof of lower bounds
%for $\dim\Hnl$, $\dim\Fnl$ and $\dim\Hnl+\Fnl$
we state formulas for values
of weight systems in~$\Hnl$ and~$\Fnl$ on certain trivalent diagrams
in Section~\ref{lemmas}.
We prove these formulas
by making calculations in the Brauer algebra~$\Brk$.
In Section~\ref{propproofs} we complete the proofs of 
Theorem~\ref{t:dimtab}, Corollary~\ref{c:hkjy}
and Theorem~\ref{t:alggen} by using
a module structure on the space of
primitive elements~$\P$ of~$\Ab$
over Vogel's algebra~$\Lambda$ (\cite{Vog}).

\section*{Acknowledgements}

I would like to thank C.-F.\ B\"odigheimer, C.\ Kassel, J.\ Kneissler, T. Mennel, 
H.\ R.\ Morton, and 
A.\ Stoimenov for helpful remarks and discussions that influenced an old version of
this article entitled
''The common specializations of the Homfly and Kauffman polynomials''. 
I thank the Graduiertenkolleg for
mathematics of the University of Bonn, the German Academic Exchange Service, and the
Schweizerischer Nationalfonds for financial support.

\section{Vassiliev invariants and link polynomials}\label{linkpolys}

A singular link is an immersion of a finite number of
oriented circles into $\R^3$
whose only singularities are transversal double points.
%We only consider tame singular links~$L$. That means that in a
%neighborhood of every point on~$L$, the link is diffeomorph to an
%interval or to intersecting intervals like the ones shown on the left side of
%Figure~\ref{f:hkpic}.
A singular link without double points is called a link. We consider
singular links up to orientation preserving diffeomorphisms of $\R^3$.
The equivalence classes of this equivalence relation are called
singular link types or by abuse of language simply singular links.
A link invariant is a map from link types into a set.
If~$v$ is a link invariant with values in an abelian group, then it
can be extended recursively
to an invariant of singular links by the local replacement rule
$v(L_\times)=v(L_+)-v(L_-)$ (see Figure~\ref{f:hkpic}).
A link invariant is called a {\em Vassiliev} invariant of degree~$n$ if it
vanishes on all singular links with~$n+1$ double points.
Let~$\Vnl$ be the vector space of~$\Q$-valued Vassiliev invariants
of degree~$n$ of links with~$\ell$ components.

\begin{figure}[!ht]
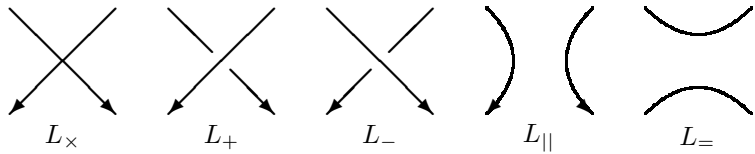

$$
\Picture{
\thicklines
\put(1,1){\vector(-1,-1){2}}
\put(-1,1){\vector(1,-1){2}}
\put(-1,-1.8){\makebox(2,0.7){$L_\times$}}
}
\qquad
\Picture{
\thicklines
\put(1,1){\vector(-1,-1){2}}
\put(-1,1){\line(1,-1){0.83}}
\put(0.17,-0.17){\vector(1,-1){0.83}}
\put(-1,-1.8){\makebox(2,0.7){$L_+$}}
}
\qquad
\Picture{
\thicklines
\put(1,1){\line(-1,-1){0.83}}
\put(-0.17,-0.17){\vector(-1,-1){0.83}}
\put(-1,1){\vector(1,-1){2}}
\put(-1,-1.8){\makebox(2,0.7){$L_-$}}
}
\qquad
\Picture{
\thicklines
\qbezier[90](1,1)(0,0)(1,-1)
\qbezier[90](-1,1)(0,0)(-1,-1)
\put(-0.9,-0.9){\vector(-1,-1){0.1}}
\put(0.9,-0.9){\vector(1,-1){0.1}}
\put(-1,-1.8){\makebox(2,0.7){$L_{\mid\mid}$}}
}
\qquad
\Picture{
\thicklines
\qbezier[90](1,1)(0,0)(-1,1)
\qbezier[90](1,-1)(0,0)(-1,-1)
\put(-1,-1.8){\makebox(2,0.7){$L_=$}}
}
$$\nopagebreak\vspace*{15pt}
\caption{Local modifications (of a diagram) of a (singular) link}\label{f:hkpic}
\end{figure}

Let us recall the definitions of the link invariants $H$,$F$,$V$,
and~$\Upsilon$ (see~\cite{HOM}, \cite{Ka2}, \cite{Jon}, and Proposition~4.7 
of~\cite{Lik}; the 
normalizations of~$H$ and~$V$ we will use are equivalent
to the original definitions).
For a link~$L$, the Homfly polynomial $H_L(x,y)\in\Z[x^{\pm 1},y^{\pm 1}]$
is given by

\begin{eqnarray}
& & x H_{L_+}(x,y)-x^{-1}H_{L_-}(x,y) =
y H_{L_{\mid\mid}}(x,y),\label{e:skein}\\
& & H_{O^k}(x,y)=\left(\frac{x-x^{-1}}{y}\right)^k.\label{e:unlink}
\end{eqnarray}

The links in
Equation~(\ref{e:skein}) are the same outside of a small ball and differ
inside this ball as shown in Figure~\ref{f:hkpic}.
The symbol~$O^k$ denotes
the trivial link with $k\geq 1$ components.

A link diagram~$L\subset \R^2$ is a generic projection of a link
together with the information which strand is the overpassing
strand at each double point of the projection.
Call a crossing of a link diagram as in~$L_+$ (see Figure~\ref{f:hkpic})
positive and a crossing as in~$L_-$ negative. Define the
{\em writhe}~$w(L)$ of a link diagram~$L$ as the number of positive crossings 
minus the number of negative crossings.
Similar to the Homfly polynomial, 
the Dubrovnik version of the {\em Kauffman} polynomial 
$F_L(x,y)\in\Z[x^{\pm 1},y^{\pm 1}]$ of a link diagram~$L$
is given by

\begin{eqnarray}
\kern-20pt xF_{L_+}(x,y)-x^{-1}F_{L_-}(x,y) \kern-8pt&=&\kern-8pt
y \left(F_{L_{\mid\mid}}(x,y)-x^{w({L_{=,or}})-w({L_{\mid\mid}})}
F_{L_{=,or}}(x,y)\right),\kern-10pt\label{e:skein2}\\
F_{O^k}(x,y) &=&\left(\frac{x-x^{-1}+y}{y}\right)^k.
\end{eqnarray}

Here the link diagrams 
$L_+,L_-,L_{\vert\vert},L_=$  
differ inside of a disk as shown in Figure~\ref{f:hkpic}
and coincide on the outside of this disk,
and~$L_{=,or}$ is the link diagram~$L_=$ %(see Figure~\ref{f:hkpic})
equipped with an arbitrary orientation of the components of the corresponding link.
The symbol~$O^k$ denotes an arbitrary diagram of 
the trivial link with~$k\geq 1$ components.
The Homfly and the Kauffman polynomials are invariants of links.
%\footnote{One can also 
%define the Kauffman polynomial first for framed links
%without orientation and then to normalize it to an invariant of oriented
%links without framing (compare~\cite{Kau}).}

Let~$\vert L\vert$ denote the number of components of a link~$L$.
For the links in Equation~(\ref{e:skein}) we
have~$\vert L_+\vert=\vert L_-\vert=\vert L_{\vert\vert}\vert\pm 1$.
Since Equations~(\ref{e:skein}) and~(\ref{e:unlink}) are sufficient to calculate~$H$ 
this implies $H_L(x,y)=(-1)^{\vert L\vert}H_L(x,-y)$
for every link~$L$.
The {\em Jones} polynomial~$V$ 
can be expressed in terms of the
Homfly polynomial as

$$
V_L(x):=H_L\left(x^2,x^{-1}-x\right)=
(-1)^{\vert L\vert}H_L\left(x^2,x-x^{-1}\right)\in\Z[x^{\pm 1}].
$$

It is easy to see that for every link~$L$ we have

\begin{equation}\label{HFval}
\widetilde{H}_L(x,y):=y^{\vert L\vert}H_L(x,y)\in \Z[x^{\pm 1},y]\ \mbox{,}
\  \widetilde{F}_L(x,y):=y^{\vert L\vert}F_L(x,y)\in \Z[x^{\pm 1},y].
\end{equation}

The link invariant~$\Upsilon$ is defined by

$$
\Upsilon_L(x)=\widetilde{H}_L(x,0)\in\Z[x^{\pm 1}].
$$

%The polynomial $\Upsilon$ cannot be obtained directly from $H$ by
%a substitution of parameters. Thus the invariants that can be obtained
%by a substitution of parameters from $H$ or $F$ depend heavily on the
%normalization used for these invariants.
%The link invariants $H$,$F$,$V$, and~$\Upsilon$ are
%{\em not} Vassiliev invariants, but
%they can be developped into
%formal power-series whose coefficients are Vassiliev
%invariants.
 
After %the following 
substitutions of parameters 
we can express~$H$ and~$F$ as

\begin{eqnarray}
H_L\left(e^{ch/2},e^{h/2}-e^{-h/2}\right) & = &
\sum_{j=0}^\infty\sum_{i=1}^{j+\vert L\vert}
p_{i,j}^{\vert L\vert}(L) c^i h^j\in\Q[c][[h]],\label{Hnpij}\\
F_L\left(e^{(c-1)h/2},e^{h/2}-e^{-h/2}\right) & = &
\sum_{j=0}^\infty
\sum_{i=1}^{j+\vert L\vert} q_{i,j}^{\vert L\vert}(L)
c^i h^j\in\Q[c][[h]],\label{Fnqij}
\end{eqnarray}

%where $c$ and $h$ are commuting variables. 
for the following reasons: Equation~(\ref{HFval}) implies
that the sum over~$i$ is limited by $j+\vert L\vert$ in these expressions and 
one sees
that no negative powers in~$h$ appear and that the sum over~$i$ starts with~$i=1$ 
directly by using the defining equations of~$H$ and~$F$ with the new parameters.
For~$j=0$ we have $p_{i,0}^{\vert L\vert}=q_{i,0}^{\vert L\vert}=\delta_{i,\vert L\vert}$, where~$\delta_{i,j}$ is~$1$ for~$i=j$ and is~$0$ otherwise.
It follows from Equations~(\ref{e:skein}) and~(\ref{e:skein2}) 
that the 
link invariants~$p_{i,n}^\ell$ and~$q_{i,n}^\ell$ are in~$\Vnl$.
%\footnote{The reason why we replaced~$x$ by the power series~$e^{ch/2}$ 
%(resp.\ $e^{(c-1)h/2}$) in Equations~(\ref{Hnpij}) and~(\ref{Fnqij}) and not by any other
%power series in~$h$ with constant term~$1$ will be given in Proposition~\ref{WglWso}.
%}. 
Define

\begin{eqnarray}
& \Hvnl=\vspan\{p_{1,n}^\ell,p_{2,n}^\ell,\ldots,p_{n+\ell,n}^\ell\}
\subseteq\Vnl, &\\
& \Fvnl=\vspan\{q_{1,n}^\ell, q_{2,n}^\ell,\ldots,
q_{n+\ell,n}^\ell\}\subseteq\Vnl. &
\end{eqnarray}

Define the invariants $y_n^\ell$, $r_n^\ell$ of links with $\ell$
components by

\begin{eqnarray}
\Upsilon_L\left(e^{h/2}\right) & = & \sum_{n=0}^\infty y_n^{\vert L\vert}(L)
h^n\in\Q[[h]],\label{Yynl}\\
V_L\left(e^{h/2}\right) & = & \sum_{n=0}^\infty r_n^{\vert L\vert}(L)
h^n\in\Q[[h]].\label{Vrnl}
\end{eqnarray}

In the following proposition we state the consequences of 
Propositions~4.7, 4.2, 4.5 of~\cite{Lik} for the versions 
of the Homfly and Kauffman polynomials from Equations~(\ref{Hnpij}) and~(\ref{Fnqij}).

\begin{prop}\label{p:hkpolys}
For all $n\geq 0,\ell\geq 1$ we have \nopagebreak{}
%where the link invariants $y_n^\ell,p_{i,j}^\ell,q_{i,j}^\ell$ are defined
%by Equations~(\ref{Yynl}), (\ref{Hnpij}), (\ref{Fnqij}), respectively.

\begin{eqnarray*}
& (1) & y_n^\ell=p_{n+\ell,n}^\ell=q_{n+\ell,n}^\ell,\\
& (2) & r_n^\ell=(-1)^\ell\sum_{i=1}^{n+\ell} 2^i p_{i,n}^\ell=
(-1/2)^n
\sum_{i=1}^{n+\ell} (-2)^i q_{i,n}^\ell,\\
& (3) & (-2)^n\sum_{i=1}^{n+\ell}4^iq_{i,n}^\ell=\sum_
{m=0}^n
\sum_{i=1}^{m+\ell}\sum_{j=1}^{n-m+\ell}(-2)^{i+j}q_{i,m}^\ell q_{j,n-m}^\ell.
\end{eqnarray*}
\end{prop}
\begin{proof}[Sketch of Proof]
(1) The following formulas for $\Upsilon$ can directly be derived
from its definition:

\begin{eqnarray*}
(a) &
x \Upsilon_{L_+}(x)-x^{-1} \Upsilon_{L_-}(x)=
\Upsilon_{L_{\mid\mid}}(x)
& \mbox{if $\vert L_+\vert < \vert L_{\mid\mid}\vert$,}\\
(b) & x \Upsilon_{L_+}(x)=x^{-1} \Upsilon_{L_-}(x) &
\mbox{if $\vert L_+\vert > \vert L_{\mid\mid}\vert$,}\\
(c) & \Upsilon_{O^k}(x)=(x-x^{-1})^k. &
\end{eqnarray*}

These relations are sufficient to calculate $\Upsilon_L(x)$
for every link~$L$. The link invariant~$\Upsilon'_L(x):=\widetilde{F}_L(x,0)$
satisfies the same Relations~(a), (b), (c) as~$\Upsilon$,
hence we have
$\widetilde{H}_L(x,0)=\Upsilon_L(x)=\Upsilon'_L(x)=\widetilde{F}_L(x,0).$
Now the formulas

$$
\widetilde{H}_L\left(e^{h/2},0\right)=\sum_{n=0}^\infty
p_{n+\vert L\vert,n}^{\vert L\vert}(L) h^n \quad\mbox{and}\quad
\widetilde{F}_L\left(e^{h/2},0\right)=\sum_{n=0}^\infty
q_{n+\vert L\vert,n}^{\vert L\vert}(L) h^n
$$

imply Part~(1) of the proposition.

\smallskip

(2) Let
$<L>(A)$ be the Kauffman bracket (see \cite{Ka1}, \cite{Kau}) defined by 

$$
<
\begin{picture}(1.2,1)(-0.6,-0.167)
\setlength{\unitlength}{8pt}
\thicklines
\put(1,1){\line(-1,-1){2}}
\qbezier[90](1,-1)(0.58,-0.58)(0.17,-0.17)
\qbezier[90](-1,1)(-0.58,0.58)(-0.17,0.17)
%\put(-1,1){\line(1,-1){0.83}}
%\put(0.17,-0.17){\line(1,-1){0.83}}
\end{picture}>
\ = \ A
<
\begin{picture}(1.2,1)(-0.6,-0.167)
\setlength{\unitlength}{8pt}
\thicklines
\qbezier[90](1,1)(0,0)(1,-1)
\qbezier[90](-1,1)(0,0)(-1,-1)
%\put(-0.9,-0.9){\vector(-1,-1){0.1}}
%\put(0.9,-0.9){\vector(1,-1){0.1}}
\end{picture}>
\ + \ A^{-1}
<
\begin{picture}(1.2,1)(-0.6,-0.167)
\setlength{\unitlength}{8pt}
\thicklines
\qbezier[90](1,1)(0,0)(-1,1)
\qbezier[90](1,-1)(0,0)(-1,-1)
\end{picture}
>
\quad , \quad <O^k>=(-A^2-A^{-2})^k.
$$\nopagebreak\vspace*{4pt}

For a link diagram~$L$
define the link invariant~$f_L(A)$ with values in $\Z[A^2,A^{-2}]$
by $f_L(A)=(-A^3)^{-w(L)}<L>(A)$, where $w(L)$ denotes the writhe of~$L$.
Then one can show that %it was proved by Kauffman that

\begin{eqnarray*}
F_L\left(e^{-3h/2},e^{h/2}-e^{-h/2}\right) = f_L\left(-e^{-h/2}\right)
& = & f_L\left(e^{-h/2}\right)
\qquad \mbox{and}\\
(-1)^{\vert L\vert}H_L\left(e^h,e^{h/2}-e^{-h/2}\right) = V_L\left(e^{h/2}\right) & = &
f_L\left(e^{h/4}\right).
\end{eqnarray*}

This implies Part~(2) of the proposition.

\smallskip

(3) With the notation of Part~(2) of the proof we have

$$
F_L(B^3,B-B^{-1})=f_L(-A^{-1})^2=F_L(A^{-3},A-A^{-1})^2,\ \mbox{where $B=A^{-2}$.}
$$

Substituting $A=e^{h/2}$ and $B=e^{\hbar/2}$ with $\hbar=-2h$ 
and comparing with Equation~(\ref{Fnqij})
gives us Part~(3) of the proposition.
\end{proof}%$\Box$

\medskip

%\begin{coro}\label{c:eleofinters}
%For all $n\geq 0,\ell\geq 1$ 
Parts~(1) and~(2) of 
Proposition~\ref{p:hkpolys} imply that $r_n^\ell,y_n^\ell\in\Hvnl\cap\Fvnl$.
%\end{coro}
In other words, the polynomials $V$ and $\Upsilon$ are common
specializations of~$H$ and~$F$.
This was the easy part of the proofs of Theorem~\ref{t:dimtab} and Corollary~\ref{c:hkjy}.
Part~(3) of Proposition~\ref{p:hkpolys} will be used in the proof of Theorem~\ref{t:alggen}.

\section{Spaces of weight systems}\label{dimtabs}

We recall the following from \cite{BN1}. A {\em trivalent diagram}
is an unoriented graph with~$\ell\geq 1$ disjointly embedded oriented
circles such that
every connected component of this graph contains at least one
oriented circle,
every vertex has valency three, and
the vertices that do not lie on an oriented circle have a cyclic
orientation. We consider trivalent diagrams up to homeomorphisms of
graphs that respect the additional data.
The degree of a trivalent diagram is defined as half of the number of
its vertices.
An example of a diagram on two circles of degree $8$ is shown in
Figure~\ref{f:trdiag}.

\begin{figure}[!ht]
$$
\begin{picture}(5.2,1.2)(-2.6,-0.1)
\thicklines
\put(-1.5,0){\circle{2}}
\put(1.5,0){\circle{2}}
\thinlines
\put(-1.5,0.5){\line(1,0){2.134}}
\put(-1,0.866){\line(1,0){2}}
\put(0,0.5){\line(0,1){0.366}}
\put(-1.5,-1){\line(0,1){2}}
\put(-0.634,-0.5){\line(-1,0){1.772}}
\put(1.5,-1){\line(0,1){1}}
\put(1.5,0){\line(1,1){0.707}}
\put(1.5,0){\line(-1,1){0.707}}
\put(-1.5,0){\line(5,-1){3}}
\end{picture}\vspace{9pt}
$$
\caption{A trivalent diagram}\label{f:trdiag}
\end{figure}
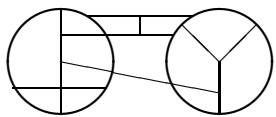

In the picture
the distinguished circles are drawn with thicker lines than the remaining part of the
diagrams.
Orientation of circles and vertices are assumed to be counterclockwise. 
Crossings in the picture do not correspond to vertices of
a trivalent diagram.
Let $\Anl$ be the $\Q$-vector space generated by trivalent diagrams of degree
$n$ on $\ell$ oriented circles together with the relations (STU), (IHX) and
(AS) shown in Figure~\ref{f:rels}.

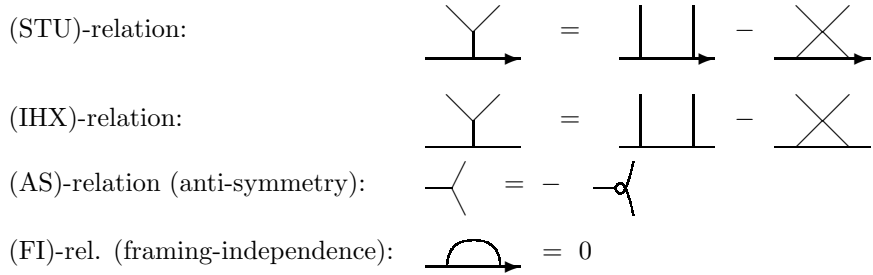
\begin{figure}[!ht]
$$
\begin{array}{ll}
\mbox{(STU)-relation:} &
\begin{picture}(2,1.1)(-1,-0.6)
\thicklines
\put(-0.9,-1){\vector(1,0){1.8}}
\thinlines
\put(0.0,-1.0){\line(0,1){0.5}}
\put(0.0,-0.5){\line(1,1){0.5}}
\put(0.0,-0.5){\line(-1,1){0.5}}
\end{picture}\quad = \quad
\begin{picture}(2,1.1)(-1,-0.6)
\thicklines
\put(-0.9,-1){\vector(1,0){1.8}}
\thinlines
\put(0.5,-1){\line(0,1){1}}
\put(-0.5,-1){\line(0,1){1}}
%\put(0,-1.4){\makebox(0,0){${}_{D_{\mid\mid}}$}}
\end{picture} \  - \ 
\begin{picture}(2,1.1)(-1,-0.6)
\thicklines
\put(-0.9,-1){\vector(1,0){1.8}}
\thinlines
\put(0.5,-1){\line(-1,1){1}}
\put(-0.5,-1){\line(1,1){1}}
\end{picture}\vspace{8pt}\\
\mbox{(IHX)-relation:} &
\begin{picture}(2,1.1)(-1,-0.6)
\thinlines
\put(-0.9,-1){\line(1,0){1.8}}
\put(0.0,-1.0){\line(0,1){0.5}}
\put(0.0,-0.5){\line(1,1){0.5}}
\put(0.0,-0.5){\line(-1,1){0.5}}
\end{picture}\quad = \quad
\begin{picture}(2,1.1)(-1,-0.6)
\thinlines
\put(-0.9,-1){\line(1,0){1.8}}
\put(0.5,-1){\line(0,1){1}}
\put(-0.5,-1){\line(0,1){1}}
%\put(0,-1.4){\makebox(0,0){${}_{D_{\mid\mid}}$}}
\end{picture} \  - \ 
\begin{picture}(2,1.1)(-1,-0.6)
\thinlines
\put(-0.9,-1){\line(1,0){1.8}}
\put(0.5,-1){\line(-1,1){1}}
\put(-0.5,-1){\line(1,1){1}}
\end{picture}\vspace{8pt}\\
\mbox{(AS)-relation (anti-symmetry):} &
\begin{picture}(0.95,0.65)(-0.35,-0.05)
\thinlines
\put(0.25,0){\line(-1,0){0.5}}
\put(0.25,0){\line(1,2){0.25}}
\put(0.25,0){\line(1,-2){0.25}}
\end{picture}\quad = \ -\quad
\begin{picture}(0.95,0.65)(-0.35,-0.05)
\thinlines
\put(0.15,0){\line(-1,0){0.4}}
\qbezier[60](0.15,0)(0.25,-0.2)(0.35,0)
\qbezier[60](0.35,0)(0.4,0.1)(0.5,0.5)
\qbezier[60](0.15,0)(0.25,0.2)(0.35,0)
\qbezier[60](0.35,0)(0.4,-0.1)(0.5,-0.5)
\end{picture}\vspace{8pt}\\
\mbox{(FI)-rel. (framing-independence):} &
\begin{picture}(2,0.7)(-1,0.15)
\thicklines
\put(-0.9,0){\vector(1,0){1.8}}
\thinlines
\qbezier[80](-0.5,0)(-0.5,0.5)(0,0.5)
\qbezier[80](0.5,0)(0.5,0.5)(0,0.5)
\end{picture}\ = \ 0
\end{array}\vspace{2pt}
$$
\caption{(STU), (IHX), (AS) and (FI)-relation}\label{f:rels}
\end{figure}

The diagrams in a relation are assumed to coincide everywhere
except for the parts we have shown.
Let~$\Abnl$ be the quotient of~$\Anl$ by the relation~(FI), also shown
in Figure~\ref{f:rels}.
A {\em weight system} is a linear map from~$\Abnl$ to a~$\Q$-vector space.

A {\em chord diagram}
is a trivalent diagram where every trivalent vertex lies
on an oriented circle. It is easy to see that $\Anl$ is spanned by chord
diagrams.
If $D$ is a chord diagram of degree $n$
on $\ell$ oriented circles, then one can construct
a singular link $L_D$ with $\ell$ components such that the preimages of
double points of $L_D$ correspond to the points of $D$ connected by a
chord.
The singular link $L_D$ described above
is not uniquely determined by $D$, but, if
$v\in\Vnl$, then the linear map
$W(v):\Abnl\longrightarrow\Q$ which sends $D$ to $v(L_D)$ is well-defined.
This defines a linear map $W:\Vnl
\longrightarrow\Hom(\Abnl,\Q)=\Abnl^*$.
Let us define the spaces

$$
\Hnl=W(\Hvnl) \qquad \mbox{and} \qquad
\Fnl=W(\Fvnl)\subseteq\Abnl^*.
$$

If $v_1\in\V_{n,\ell}$ and $v_2\in\V_{m,\ell}$, then the link invariant
$v_1 v_2$
defined by $(v_1 v_2)(L)=v_1(L)v_2(L)$ is in $\V_{n+m,\ell}$.
Weight systems are multiplied by using the algebra structure dual to the
coalgebra structure of $\bigoplus_{n=0}^\infty \Abnl$ (see~\cite{BN1}).
The following proposition is a well-known consequence of a theorem of Kontsevich 
(see Proposition~2.9 of
\cite{BNG} and Theorem~7.2 of~\cite{KaT}, Theorem~10 of~\cite{LM3} or
\cite{LM1}, \cite{LM2}).

\begin{prop}\label{WglWsoIntegral}
For all $\ell\geq 1$ there exists an isomorphism of algebras\nopagebreak

$$
Z^*:\bigoplus_{n=0}^\infty\Abnl^*\longrightarrow
\bigcup_{n=0}^\infty \Vnl
$$\nopagebreak

such that for all $n\geq 0$ we have\nopagebreak

$$
Z^*\circ W_{\vert (\Hvnl+\Fvnl)}=\id_{(\Hvnl+\Fvnl)}.
$$

%sending
%$W(p_{i,j}^\ell)$ back to $p_{i,j}^\ell$ (resp.
%$W(q_{i,j}^\ell)$
%to $q_{i,j}^\ell$).
\end{prop}

This proposition reduces the study of $\Hvnl$ and
$\Fvnl$ to that of $\Hnl$ and $\Fnl$: we have the following corollary.

\begin{coro}\label{c:dimtab}
For all $n\geq 0$ and $\ell\geq 1$ we have 

\begin{eqnarray*}
\dim\Hnl & = & \dim\Hvnl,\\
\dim\Fnl & = & \dim\Fvnl,\\ 
\dim(\Hnl+\Fnl) & = & \dim(\Hvnl+\Fvnl).
\end{eqnarray*}
\end{coro}

We will often use Corollary~\ref{c:dimtab} without referring to it.

\section{Upper bounds for $\dim\Hnl$ and $\dim\Fnl$}\label{upperb}

Let us recall the explicit descriptions of $W(p_{i,j}^\ell)$ and
$W(q_{i,j}^\ell)$ from~\cite{BN1}.
Let $D$ be a trivalent diagram. Cut it into pieces along small circles
around each vertex. Then replace the simple parts as shown in
Figure~\ref{f:Wgl}.

\begin{figure}[!ht]
$$
\begin{picture}(1,1)(-0.5,-0.1)
\thicklines
\put(-0.4,1){\vector(0,-1){2}}
\thinlines
\put(-0.4,0){\line(1,0){0.8}}
\end{picture}\ \leadsto \
\begin{picture}(1,1)(-0.5,-0.1)
\thicklines
\put(-0.4,1){\line(0,-1){0.8}}
\put(-0.4,0.2){\line(1,0){0.8}}
\put(-0.4,-0.2){\line(1,0){0.8}}
\put(-0.4,-0.2){\line(0,-1){0.8}}
\end{picture}\qquad\quad
\begin{picture}(0.3,1)(-0.15,-0.1)
\thinlines
\put(0,1){\line(0,-1){2}}
\end{picture}\ \leadsto \
\begin{picture}(1,1)(-0.5,-0.1)
\thicklines
\put(-0.2,1){\line(0,-1){2}}
\put(0.2,1){\line(0,-1){2}}
\end{picture}\qquad\quad
\begin{picture}(1.7,1)(-0.85,-0.1)
\thinlines
\put(-0.75,0){\line(1,0){1}}
\put(0.25,0){\line(1,2){0.5}}
\put(0.25,0){\line(1,-2){0.5}}
\end{picture}\ \leadsto \
\begin{picture}(1.7,1)(-0.85,-0.1)
\thicklines
\put(-0.75,0.2){\line(1,0){0.8}}
\put(-0.75,-0.2){\line(1,0){0.8}}
\put(0.05,0.2){\line(1,2){0.4}}
\put(0.45,0){\line(1,2){0.5}}
\put(0.05,-0.2){\line(1,-2){0.4}}
\put(0.45,0){\line(1,-2){0.5}}
\end{picture} -
\begin{picture}(1.7,1)(-0.85,-0.1)
\thicklines
\put(-0.75,-0.2){\line(3,2){1.8}}
\put(-0.75,0.2){\line(3,-2){1.8}}
\put(0.55,1){\line(0,-1){2}}
\end{picture}
$$\vspace{8pt}
\caption{The map $W_\gl$}\label{f:Wgl}
\end{figure}

Glue the substituted parts together. Sums of parts of diagrams are glued
together after multilinear expansion.
The result is a linear combination of unions of circles. Replace each circle
by a formal parameter $c$ and call the resulting polynomial $W_\gl(D)$. It
is well-known that this procedure determines a linear map
$W_\gl:\Anl\longrightarrow \Q[c]$ (see \cite{BN1}, Exercise~6.36).
Proceeding with
the replacement patterns
shown in Figure~\ref{f:Wso}, we get
the linear map $W_\so:\Anl\longrightarrow\Q[c]$.

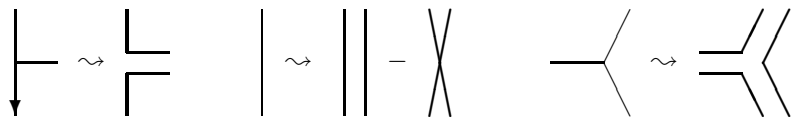
\begin{figure}[!ht]
$$
\begin{picture}(1,1)(-0.5,-0.1)
\thicklines
\put(-0.4,1){\vector(0,-1){2}}
\thinlines
\put(-0.4,0){\line(1,0){0.8}}
\end{picture}\ \leadsto \
\begin{picture}(1,1)(-0.5,-0.1)
\thicklines
\put(-0.4,1){\line(0,-1){0.8}}
\put(-0.4,0.2){\line(1,0){0.8}}
\put(-0.4,-0.2){\line(1,0){0.8}}
\put(-0.4,-0.2){\line(0,-1){0.8}}
\end{picture}\qquad\quad
\begin{picture}(0.3,1)(-0.15,-0.1)
\thinlines
\put(0,1){\line(0,-1){2}}
\end{picture}\ \leadsto \
\begin{picture}(1,1)(-0.5,-0.1)
\thicklines
\put(-0.2,1){\line(0,-1){2}}
\put(0.2,1){\line(0,-1){2}}
\end{picture} -
\begin{picture}(1,1)(-0.5,-0.1)
\thicklines
\put(-0.2,1){\line(1,-5){0.4}}
\put(0.2,1){\line(-1,-5){0.4}}
\end{picture}\qquad\quad
\begin{picture}(1.7,1)(-0.85,-0.1)
\thinlines
\put(-0.75,0){\line(1,0){1}}
\put(0.25,0){\line(1,2){0.5}}
\put(0.25,0){\line(1,-2){0.5}}
\end{picture}\ \leadsto \
\begin{picture}(1.7,1)(-0.85,-0.1)
\thicklines
\put(-0.75,0.2){\line(1,0){0.8}}
\put(-0.75,-0.2){\line(1,0){0.8}}
\put(0.05,0.2){\line(1,2){0.4}}
\put(0.45,0){\line(1,2){0.5}}
\put(0.05,-0.2){\line(1,-2){0.4}}
\put(0.45,0){\line(1,-2){0.5}}
\end{picture}
$$\vspace{8pt}
\caption{The map $W_\so$}\label{f:Wso}
\end{figure}

For a trivalent diagram $D$, define the linear combination of
trivalent diagrams~$\iota(D)$ by
replacing each chord as
shown in Figure~\ref{f:iota}. Connected components of~$D\setminus{S^1}^{\amalg\ell}$
with an internal trivalent vertex stay as they are. 

\begin{figure}[!ht]
$$
\begin{picture}(1.6,1)(-0.8,-0.1)
\thicklines
\put(-0.7,1){\vector(0,-1){2}}
\put(0.7,-1){\vector(0,1){2}}
\thinlines
\put(-0.7,0){\line(1,0){1.4}}
\end{picture}\ \leadsto \
\begin{picture}(1.6,1)(-0.8,-0.1)
\thicklines
\put(-0.7,1){\vector(0,-1){2}}
\put(0.7,-1){\vector(0,1){2}}
\thinlines
\put(-0.7,0){\line(1,0){1.4}}
\end{picture} -\frac{1}{2}\left(\
\begin{picture}(1.6,1)(-0.8,-0.1)
\thicklines
\put(-0.7,1){\vector(0,-1){2}}
\put(0.7,-1){\vector(0,1){2}}
\thinlines
\qbezier[80](-0.7,0.4)(-0.3,0.4)(-0.3,0)
\qbezier[80](-0.7,-0.4)(-0.3,-0.4)(-0.3,0)
\end{picture} +
\begin{picture}(1.6,1)(-0.8,-0.1)
\thicklines
\put(-0.7,1){\vector(0,-1){2}}
\put(0.7,-1){\vector(0,1){2}}
\thinlines
\qbezier[80](0.7,0.4)(0.3,0.4)(0.3,0)
\qbezier[80](0.7,-0.4)(0.3,-0.4)(0.3,0)
\end{picture}\ \right)
$$\vspace{8pt}
\caption{The deframing map $\iota$}\label{f:iota}
\end{figure}
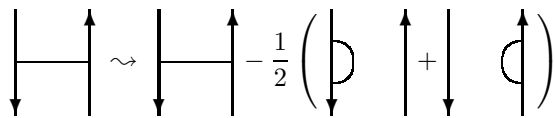

This definition determines a linear map
$\iota:\Abnl\longrightarrow \Anl$, such that $\pi\circ\iota=\id$ where
$\pi:\Anl\longrightarrow \Abnl$ denotes the canonical projection (compare
\cite{BN1}, Exercise~3.16).
By the following proposition (\cite{BN1}, Chapter~6.3) 
the weight systems
$\Wgl=W_\gl\circ\iota$ and~$\Wso=W_\so\circ\iota$ belong to the Homfly and
Kauffman polynomials.

\begin{prop}\label{WglWso}
For all $n\geq 0,i,\ell\geq 1$ the weight system
$W(p_{i,n}^\ell)$ (resp.\ $W(q_{i,n}^\ell)$) is equal to the
coefficient of $c^i$ in $\Wgl_{\vert\Abnl}$ (resp.\ $\Wso_{\vert\Abnl}$).
\end{prop}

The direct description of $W(p_{i,n}^\ell)$ and $W(q_{i,n}^\ell)$ from the
proposition above
will simplify the computation of dimensions.

\begin{lemma}\label{l:upperbo}
(1) For all $n,\ell\geq 1$ we have\nopagebreak

%\begin{equation}\label{dimhnl}
$$
\dim \Hnl\leq
\left\{
\begin{array}{ll}
n & \mbox{if $n<\ell$,}\\
\left[\frac{n-1+\ell}{2}\right] & \mbox{if $n\geq \ell$.}
\end{array}
\right.
%\end{equation}
$$

(2) For all $n,\ell\geq 1$ we have

%\begin{equation}\label{dimknl}
$$
\dim \Knl\leq\left\{
\begin{array}{ll}
n-1 & \mbox{if $\ell=1$,}\\
2n-1 & \mbox{if $\ell\geq 2$ and $n\leq \ell$,}\\
n+\ell-1 & \mbox{if $\ell\geq 2$ and $n\geq \ell$}.
\end{array}
\right.
%\end{equation}
$$

%(3) For all $n\geq 3$ we have $\dim \Fn1_{\vert\P_n}\leq n-2$.
\end{lemma}
\begin{proof}
In the proof
$D$ will denote a chord diagram of degree $n\geq 1$ on $\ell$ circles.

(1) If $n\geq \ell$, then we get $[(n-1+\ell)/2]$
as an upper bound for $\dim \Hnl$
by the following observations:

(a) The polynomial $\Wgl(D)$ has 
degree $\leq n+\ell$ and
vanishing constant term 
because the number of circles can at most increase by
one with each replacement of a chord as shown in Figure~\ref{f:Wgl}, and
there remains always at least one circle.

(b) The coefficients of $c^{n+\ell-1-2i}$ ($i=0,1,\ldots$) vanish because
the number of circles changes by $\pm 1$ with each
replacement of a chord as shown in Figure~\ref{f:Wgl}.

(c) We have $\Wgl(D)(1)=0$ because $W_\gl(D')(1)=1$ for each chord diagram
$D'$ and $\iota(D)$ is a linear combination of chord diagrams $D'$ having
$0$ as sum of their coefficients.

\smallskip

If $D$ is a chord diagram of degree $n<\ell$, then by similar arguments
$\Wgl(D)$ is a linear combination of $c^{\ell-n},c^{\ell-n+2},
\ldots,c^{\ell+n}$
with $\Wgl(D)(1)=0$. This implies the upper bound for $\dim\Hnl$.

\smallskip

(2) If $n\geq \ell$, then by the same arguments as above $\Wso(D)$ is a
polynomial of degree $\leq n+\ell$ with vanishing constant term and
$\Wso(D)(1)=0$. This implies $\dim\Fnl\leq n+\ell-1$ in this case.

If $\ell=1$, then for chord diagrams $D'$ of degree $n$ the value
$W_\so(D')(2)$ is constant because $\so_2$ is an abelian
Lie algebra (see \cite{BN1}). This implies $\Wso(D)(2)=0$ and hence
$\dim \Fnl\leq n-1$ in this case.

If $\ell\geq 2$ and $n<\ell$, then the coefficient of~$c^{\ell-n}$
in~$\Wso(D)$
is~$0$ by the following argument: 
Assume that a chord diagram~$D'$ has the minimal possible number of~$\ell-n$ 
connected components (in other
words, if we contract the oriented
circles of~$D'$ to points, then the resulting graph
is a forest). Then we see that $W_\so(D')=0$ by using Figure~\ref{f:Wso}.
Hence~$\Wso(D)$ is a linear combination of
$c^{\ell-n+1},c^{\ell-n+2},\ldots,c^{\ell+n}$ with $\Wso(D)(1)=0$.
This completes the proof of the upper bounds for $\dim\Fnl$.
\end{proof}%$\Box$

\section{The Brauer algebra and values of $\Wgl$ and $\Wso$}\label{lemmas}

In order to find lower bounds for $\dim \Hnl$, $\dim \Fnl$ and
$\dim (\Hnl+\Fnl)$, we shall evaluate the weight systems
$\Wgl$ and $\Wso$ on sufficiently many trivalent diagrams.
Let $\omega_k,L_k,C_k,T_k$ be the diagrams of degree $k$
shown in Figure~\ref{f:oLCT}.
%\footnote{We have chosen the diagrams $\omega_k$
%for our computations because they
%are the simplest diagrams~$D$ representing primitive elements such that
%$D\setminus S^1$ is not a tree.
%We have chosen the diagrams $L_k$ and $C_k$ because for
%calculations in low degrees $n$ we have to connect as many circles as
%possible.}
%The choice of $T_k$ is useful
%because the calculations will show that we need each degree $k$ and
%for all number $l\geq 2$
%of oriented circles a diagram $D$ such that
%$\iota(D)$ is not a linear combination of trivalent diagrams 
%with an internal trivalent vertex. The disjoint union
%of $l-2$ oriented circles with $T_k$ turns out to have this property.}

\begin{figure}[!ht]
$$
\begin{array}{cl}
%\begin{eqnarray*}
\omega_k = 
\Picture{
\DottedCircle
\Arc[0]\Arc[1]\Arc[2]\Arc[3]\Arc[4]\Arc[5]
\Arc[6]\Arc[7]\Arc[8]\Arc[10]\Arc[11]
\Endpoint[7]\Endpoint[8]\Endpoint[10]\Endpoint[11]
\thinlines
\put(-0.866,-0.5){\line(1,2){0.25}}
\put(-0.5,-0.866){\line(1,3){0.289}}
\put(0.866,-0.5){\line(-1,2){0.25}}
\put(0.5,-0.866){\line(-1,3){0.289}}
\put(-1,-1){\makebox(2,1){$\scriptstyle\ldots$}}
\put(0,0.2){\oval(1.5,0.4)}
} &
L_k = 
\begin{picture}(8,2)(-4,-0.1)
\thicklines
\put(-3.4,0){\circle{1}}
\put(-1.7,0){\circle{1}}
\put(3.4,0){\circle{1}}
\put(1.7,0){\circle{1}}
\thinlines
\put(-2.9,0){\line(1,0){0.7}}
\put(-1.2,0){\line(1,0){0.7}}
\put(2.9,0){\line(-1,0){0.7}}
\put(1.2,0){\line(-1,0){0.7}}
\put(-1,-0.5){\makebox(2,1){$\cdots$}}
\end{picture}\\[5pt]
T_k =
\begin{picture}(5.2,2.2)(-2.6,-0.1)
\thicklines
\put(-1.5,0){\circle{2}}
\put(1.5,0){\circle{2}}
\thinlines
\put(-0.634,0.5){\line(1,0){1.268}}
\put(-1,0.866){\line(1,0){2}}
\put(-0.634,-0.5){\line(1,0){1.268}}
\put(-1,-0.866){\line(1,0){2}}
\put(-1,-0.4){\makebox(2,1){$\vdots$}}
\end{picture} &
C_k = 
\begin{picture}(8,1.5)(-4,-0.1)
\thicklines
\put(-3.4,0){\circle{1}}
\put(-1.7,0){\circle{1}}
\put(1.7,0){\circle{1}}
\put(3.4,0){\circle{1}}
\thinlines
\put(-2.9,0){\line(1,0){0.7}}
\put(-1.2,0){\line(1,0){0.7}}
\put(1.2,0){\line(-1,0){0.7}}
\put(2.9,0){\line(-1,0){0.7}}
\qbezier[100](-3.4,0.5)(-3.4,0.8)(-2,0.8)
\qbezier[100](3.4,0.5)(3.4,0.8)(2,0.8)
\put(-2,0.8){\line(1,0){4}}
\put(-1,-0.5){\makebox(2,1){$\cdots$}}
\end{picture}\\[8pt]
%\end{eqnarray*}\vspace{8pt}
\end{array}
$$
\caption{The diagrams $\omega_k$, $L_k$, $C_k$, $T_k$}\label{f:oLCT}
\end{figure}
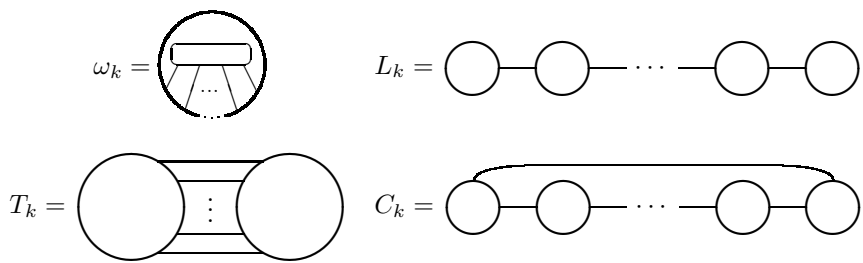

For technical reasons we extend this definition by setting $L_0=C_0=T_0=S^1$ and $C_1=L_1$.
An important ingredient in the proofs of Theorems~\ref{t:dimtab} and~\ref{t:alggen} 
is the following lemma.

\begin{lemma}\label{valofdiags}
(1) For all $k\geq 2$ we have\nopagebreak

\begin{eqnarray*}
\Wgl(\omega_k) & = & \left\{
\begin{array}{ll}
c^{k+1}+c^3-2c & \mbox{if $k$ is even,}\\
c^{k+1}-c^2 & \mbox{if $k$ is odd, and}
\end{array}
\right.\\ 
\Wso(\omega_k) & = & c(c-1)(c-2)R_k(c),
\end{eqnarray*}\nopagebreak

where $R_k$ is a polynomial with $R_k(0)\not=0$. If $k=2$, then
$R_2= 2$, and if $k\not=3$, then $R_k(2)\not=0$.

\smallskip

(2) For all $k\geq 1$ we have\nopagebreak

$$
\begin{array}{ll}
%\begin{eqnarray*}
\Wgl(L_k) = c(1-c^2)^k, &
\Wso(L_k) = c^{k+1}(1-c)^k,\\
\Wgl(T_k) = (-c)^k(c^2-1), &
\Wso(T_k) = c(c-1)Q_k(c),\\
& \Wso(C_k) = c(c-1)P_k(c),
%\end{eqnarray*}\nopagebreak
\end{array}
$$\nopagebreak

where $P_k$ and $Q_k$ are polynomials in $c$ such that for $k\geq 2$ we have
$P_k(0)\not=0$, $Q_k(0)=2^{k-1}$, and $Q_k(2)=(-2)^k$.
\end{lemma}

In the proof of the lemma we will determine
the polynomials $P_k$, $Q_k$, and $R_k$ explicitly,
which will be helpful to us for
calculations in low degrees. For the main parts of the proofs of Theorems~\ref{t:dimtab}
and~\ref{t:alggen} it will be sufficient to
know the properties of these polynomials stated in the lemma.
We do not need to know the value of~$\Wgl(C_k)$.

In the proof of Lemma~\ref{valofdiags}
we use the Brauer algebra (\cite{Bra}) on~$k$ strands~$\Brk$.
As a $\Q[c]$-module $\Brk$ has a basis in one-to-one correspondence
with involutions without fixed-points of the set
$\{1,\ldots,k\}\times\{0,1\}$.
We represent a basis element corresponding to an involution $f$ graphically
by connecting the points~$(i,j)$ and~$f(i,j)$ by a curve in~$\R\times[0,1]$.
Examples 
are the
diagrams $u_-$, $x_+$, $x_-$, $u_+=d$, $e$, $f$, $g$, $h$ in
Figures~\ref{f:fourpos} and~\ref{f:fivepos}.

\begin{figure}[!ht]
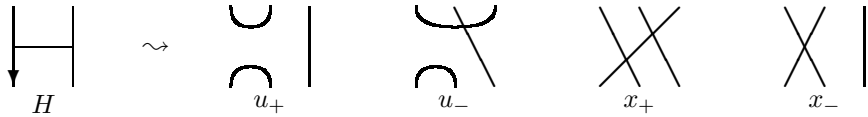

{\unitlength=15pt
$$
\Picture{
\thicklines
\put(-0.75,1){\vector(0,-1){2}}
\thinlines
\put(0.75,1){\line(0,-1){2}}
\put(-0.75,0){\line(1,0){1.5}}
\put(-1,-1.8){\makebox(2,0.7){$H$}}
} \qquad \leadsto \qquad
\Picture{
\bline[0,1]\bline[2,3]\bline[4,5]
\put(-1,-1.8){\makebox(2,0.7){$u_+$}}
}
\qquad\qquad
\Picture{
\bline[0,1]\bline[2,4]\bline[3,5]
\put(-1,-1.8){\makebox(2,0.7){$u_-$}}
}
\qquad\qquad
\Picture{
\bline[0,3]\bline[1,5]\bline[2,4]
\put(-1,-1.8){\makebox(2,0.7){$x_+$}}
}
\qquad\qquad
\Picture{
\bline[0,4]\bline[1,5]\bline[2,3]
\put(-1,-1.8){\makebox(2,0.7){$x_-$}}
}
$$\vspace*{8pt}}
\caption{Elements of $\Br3$ needed to calculate
$W_\so(\omega_k)$}\label{f:fourpos}
\end{figure}

\begin{figure}[!ht]
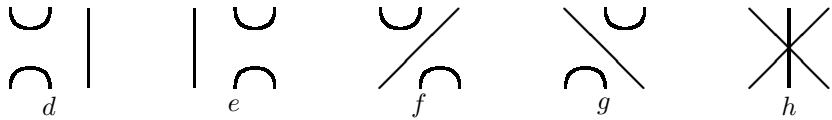

{\unitlength=15pt
$$
\Picture{
\bline[0,1]\bline[2,3]\bline[4,5]
\put(-1,-1.8){\makebox(2,0.7){$d$}}
}
\qquad\qquad
\Picture{
\bline[0,5]\bline[1,2]\bline[3,4]
\put(-1,-1.8){\makebox(2,0.7){$e$}}
}
\qquad\qquad
\Picture{
\bline[0,3]\bline[1,2]\bline[4,5]
\put(-1,-1.8){\makebox(2,0.7){$f$}}
}
\qquad\qquad
\Picture{
\bline[0,1]\bline[2,5]\bline[3,4]
\put(-1,-1.8){\makebox(2,0.7){$g$}}
}
\qquad\qquad
\Picture{
\bline[0,3]\bline[1,4]\bline[2,5]
\put(-1,-1.8){\makebox(2,0.7){$h$}}
}
$$\vspace*{8pt}}
\caption{Diagrams needed to calculate $W_\gl(\omega_k)$}\label{f:fivepos}
\end{figure}

The product of basis vectors $a$ and $b$ is
defined graphically by placing~$a$ onto the top of $b$,
by gluing the lower points $(i,0)$ of $a$ to the
upper points $(i,1)$ of $b$, and
by introducing the relation that a circle is equal to
the formal parameter $c$ of the ground ring $\Q[c]$.
We have a map $\tr:\Brk\longrightarrow \Q[c]$, called trace, that
is defined graphically 
by connecting the vertices $(i,0)$ and $(i,1)$ of a diagram
by curves, and by replacing each circle by the indeterminate~$c$.
As an example, the trace of the diagram $x_+u_-$ is
shown in Figure~\ref{f:trace}.

\begin{figure}[!ht]
{\unitlength=15pt
$$
\tr\left(
\begin{picture}(2,2)(-1,-0.167) 
\bline[0,1]\bline[2,5]\bline[3,4]
\end{picture}
\right) \ = \
\begin{picture}(4,2)(-1,-0.167) 
\btrace
\bline[0,1]\bline[2,5]\bline[3,4]
\end{picture}
\ = \ c
$$\vspace*{6pt}}
\caption{The trace of a diagram}\label{f:trace}
\end{figure}

The elements $u_+, u_-,x_+,x_-$ arise among others
when the replacement rules
belonging to $W_\so$ (see Figure~\ref{f:Wso}) are applied to the part $H$
(see Figure~\ref{f:fourpos}) of a trivalent diagram.
Similarly, the elements $d$ and $h$ arise when we apply the replacement
rules belonging to $W_\gl$ (see Figure~\ref{f:Wgl}) to the part $H$ of a
trivalent diagram.
We have $\iota(\omega_k)=\omega_k$ (see Figure~\ref{f:iota}) because the
diagram $\omega_k$ contains no chords.
The proof of the following lemma is now straightforward.

\begin{lemma}\label{brauer}
The following two formulas hold:

$$
\Wgl(\omega_k)=\tr\left((d-h)^k\right)\qquad\mbox{and}\qquad
\Wso(\omega_k)=\tr\left((u_+-u_-+x_+-x_-)^k\right).
$$
\end{lemma}

Now we can prove Lemma~\ref{valofdiags} by making
calculations in the Brauer algebra.

\medskip

\begin{proof}[Proof of Lemma~\ref{valofdiags}] 
(1) With the elements $u_\pm,x_\pm\in\Br3$ shown in Figure~\ref{f:fourpos}
we define $\Du=u_+-u_-$ and $\Dx=x_+-x_-$.
It is easy to verify that

\begin{equation}\label{multux}
(\Du+\Dx)\Du=(c-2)\Du\mbox{ and }\Dx^3=\Dx^2+2\Dx.
\end{equation}

In view of the expression for $x^3$ it is clear that $x^k$ can be
expressed as a linear combination of $x$ and $x^2$:

\begin{equation}\label{dnx}
d_k\Dx^2+e_k\Dx = \Dx^k. 
\end{equation}

It can be shown by induction that the 
sequence of pairs $(d_k,e_k)_{k\geq 1}$ is given by $(d_1,e_1)=(0,1)$
and $(d_{k+1},e_{k+1})=(d_k+e_k,2d_k)$.
We deduce

\begin{eqnarray}
& & d_1-e_1=-1, d_{k+1}-e_{k+1}=d_k+e_k-2d_k=e_k-d_k\nonumber\\
& & \Rightarrow 
d_k-e_k = (-1)^k,\label{dminuse}\\
& & d_{k+1}+(-1)^k =(d_k+e_k)+(d_k-e_k) = 2d_k. \label{dnplusone}
\end{eqnarray}

By Equations~(\ref{multux}) and~(\ref{dnx}) we have

\begin{eqnarray}\label{uxn}
& &(\Du+\Dx)^k =
\Dx^k+\sum_{i=0}^{k-1} (\Du+\Dx)^i\Du\Dx^{k-i-1} \nonumber\\
& & = d_k\Dx^2+e_k\Dx+(c-2)^{k-1}\Du
+\sum_{i=0}^{k-2} (c-2)^i(d_{k-i-1}\Du\Dx^2+e_{k-i-1}\Du\Dx). 
\end{eqnarray}

It is easy to see that

\begin{equation}\label{eqtrace}
\tr\left(\Dx^2\right)/(c-1)=-\tr(\Dx)=\tr(\Du)=-\tr(\Du\Dx)=
\tr\left(\Du\Dx^2\right)=c^2-c.
\end{equation}

Applying the trace to Equation~(\ref{uxn}) yields by
Lemma~\ref{brauer} and Equations~(\ref{dminuse}) and~(\ref{eqtrace}):

\begin{eqnarray*}
&&\Wso(\omega_k)\\
 & = & (c^2-c)\left[d_k(c-1)-e_k+(c-2)^{k-1}+
\sum_{i=0}^{k-2} (c-2)^i(d_{k-i-1}-e_{k-i-1})\right]\\
 & = & (c^2-c)\left[d_k(c-2)+(-1)^k-
\sum_{i=0}^{k-1} (-1)^{k-i}(c-2)^i\right]\\
 & = & (c^2-c)(c-2)\left[(d_k+(-1)^k)+
\sum_{i=1}^{k-2} (-1)^{k-i}(c-2)^i\right].
\end{eqnarray*}

Define the sequence $(a_k)_{k\geq 2}$ inductively by $a_2=2$ and
$a_{k+1}=2 a_k-4(-1)^k$.
We have $a_2=d_2+(-1)^2$ and by definition of $a_k$, induction and
Equation~(\ref{dnplusone}) also

$$
a_{k+1}=2 a_k-4(-1)^k=2(d_k+(-1)^k)-4(-1)^k=d_{k+1}+(-1)^{k+1}.
$$

This implies $\Wso(\omega_k)=c(c-1)(c-2)R_k(c)$ with

$$
R_k(c)=a_k+\sum_{i=1}^{k-2}(-1)^{k-i}(c-2)^i.
$$

The properties of $R_k$ stated in the lemma are satisfied because by a simple
computation we have $R_k(2)=a_k>0$ for $k\not=3$ and

$$
R_k(0)=a_k+(-1)^k(2^{k-1}-2)\equiv 2\ \mod\;4.
$$

\medskip

We only give a sketch of the proof of the formula for
$\Wgl(\omega_k)$.
Let $d,e,f,g,h$ be the elements of $\Br3$ shown
in Figure~\ref{f:fivepos}.
Then one can prove by induction on $k$ that

$$
(d-h)^{2k+1}=c^{2k}d -h +\sum_{i=0}^{k-1}c^{2i}(d+e)-c^{2i+1}(f+g).
$$

Using Lemma~\ref{brauer}
this formula allows to conclude by
distinguishing whether $k$ is even or odd.

\medskip

(2) Let $a,b,{\bf 1}$ be the elements of $\Brt$ shown in Figure~\ref{f:ab}.

\begin{figure}[!ht]
$$
a=\begin{picture}(1,0.7)(-0.5,-0.1)
\thicklines
\qbezier[60](-0.4,0.6)(-0.4,0.2)(0,0.2)
\qbezier[60](0.4,0.6)(0.4,0.2)(0,0.2)
\qbezier[60](-0.4,-0.6)(-0.4,-0.2)(0,-0.2)
\qbezier[60](0.4,-0.6)(0.4,-0.2)(0,-0.2)
\end{picture}\qquad\quad
b=\begin{picture}(1,0.7)(-0.5,-0.1)
\thicklines
\put(-0.4,0.6){\line(2,-3){0.8}}
\put(0.4,0.6){\line(-2,-3){0.8}}
\end{picture}\qquad\quad
{\bf 1}=\begin{picture}(1,0.7)(-0.5,-0.1)
\thicklines
\put(-0.4,0.6){\line(0,-1){1.2}}
\put(0.4,0.6){\line(0,-1){1.2}}
\end{picture}\vspace{6pt}
$$
\caption{Diagrams in $\Brt$}\label{f:ab}
\end{figure}

Then we have $ab=ba=a$, $a^2=ca$, $b^2={\bf 1}$,
$\tr(a)=\tr(b)=c$, $\tr({\bf 1})=c^2$, and by convention $(a-b)^0={\bf 1}$.
This implies for $k\geq 1$ that

\begin{eqnarray*}
\Wso(T_k) & = & \tr\left((a-b+{\bf 1}-c{\bf 1})^k\right)\\
& = & \tr\left[\sum_{i=0}^k \bino(k,i) (1-c)^{k-i}(a-b)^i\right]\\
& = & \tr\left\{\sum_{i=0}^k \bino(k,i)(1-c)^{k-i}\left[ (-b)^i+\sum_{j=1}^i\bino(i,j)
c^{j-1}(-1)^{i-j}a\right]\right\}\\
& = & \sum_{i=0}^k \bino(k,i)(1-c)^{k-i}\left[
\tr\left((-b)^i\right) +(c-1)^i -(-1)^i \right]\\
& = & \sum_{i=0}^k \bino(k,i)(1-c)^{k-i}\left[
\tr\left((-b)^i\right) -(-1)^i \right]\\
& = & %\sum_{\mbox{\tiny $i=0$} \atop
%\mbox{\tiny $i$ even}}^k
\sum_{\substack{0\leq i\leq k \\ \mbox{\tiny $i$ even}}}
\bino(k,i)(1-c)^{k-i}(c^2-1)
+%\sum_{\mbox{\tiny $i=1$}\atop\mbox{\tiny $i$ odd}}^k
\sum_{\substack{1\leq i\leq k \\ \mbox{\tiny $i$ odd}}}
\bino(k,i)(1-c)^{k-i}(1-c)\\
& = & \sum_{i=0}^k \bino(k,i) (1-c)^{k-i}(c-1)(-1)^i+
c%\sum_{\mbox{\tiny $i=0$}\atop
%\mbox{\tiny $i$ even}}^k 
\sum_{\substack{0\leq i\leq k \\ \mbox{\tiny $i$ even}}}
\bino(k,i)
(1-c)^{k-i}(c-1)\\
& = & c (c-1) \left[ -(-c)^{k-1}+
%\sum_{\mbox{\tiny $i=0$}\atop
%\mbox{\tiny $i$ even}}^k
\sum_{\substack{0\leq i\leq k \\ \mbox{\tiny $i$ even}}}
\bino(k,i)
(1-c)^{k-i}\right].
\end{eqnarray*}

Now one checks the properties of $Q_k$ using the last expression for
$\Wso(T_k)$.
The remaining formulas follow by easy computations.
For example, $\Wso(L_k)$ is given by the value of~$W_\so$ on the diagrams in
$\iota(L_k)$ where no chord connects two different circles.
Furthermore, one can show for~$k\geq 2$ that

$$
\Wso(C_k)=\tr\left(({\bf 1}-b)^k\right)+(1-c)\Wso(L_{k-1})=
c(c-1)\left(2^{k-1}-c^{k-1}(1-c)^{k-1}\right)\kern-2pt.
$$

The property $P_k(0)\not=0$ from the lemma is obvious from the formula above.
\end{proof}%$\Box$

\section{Completion of proofs using Vogel's algebra}\label{propproofs}

In the case of diagrams on one oriented circle, 
the coalgebra structure of $\Ab=\bigoplus_{n=0}^\infty \Ab_{n}$
can be extended to a Hopf algebra structure (see~\cite{BN1}).
The primitive elements~$\P$ of~$\Ab$ are spanned by diagrams~$D$
such that~$D\setminus S^1$ is connected, where~$S^1$ denotes the oriented
circle of $D$.
Vogel defined an algebra $\Lambda$
which acts on primitive elements (see \cite{Vog}).
The diagrams
$t$ and $x_3$ shown in Figure~\ref{f:tx3x5} represent
elements of~$\Lambda$.
%(we assume that the univalent vertices of $t$, $x_3$,
%$x_5$ are labeled by $\{1,2,3\}$ as in \cite{Vog}).

\begin{figure}[!ht]
$$
t=\begin{picture}(2.2,0.6)(-1,-0.1)
\put(-0.9,0.5){\line(1,0){1}}
\put(-0.9,-0.5){\line(1,0){1}}
\put(-0.4,-0.5){\line(0,1){1}}
\qbezier[80](0.1,0.5)(0.6,0.5)(0.6,0)
\qbezier[80](0.1,-0.5)(0.6,-0.5)(0.6,0)
\put(0.6,0){\line(1,0){0.5}}
\end{picture}\qquad\quad
x_3=\begin{picture}(2.2,0.6)(-1,-0.1)
\put(-0.9,0.5){\line(1,0){1}}
\put(-0.9,-0.5){\line(1,0){1}}
\put(-0.4,-0.5){\line(0,1){1}}
\qbezier[80](0.1,0.5)(0.6,0.5)(0.6,0)
\qbezier[80](0.1,-0.5)(0.6,-0.5)(0.6,0)
\put(0.6,0){\line(1,0){0.5}}
\put(0.1,-0.5){\line(0,1){1}}
\put(-0.4,0){\line(1,0){0.5}}
\end{picture}%\qquad\quad
%x_5=\begin{picture}(2.2,0.6)(-1,-0.1)
%\put(-0.9,0.5){\line(1,0){1}}
%\put(-0.9,-0.5){\line(1,0){1}}
%\put(-0.4,-0.5){\line(0,1){1}}
%\qbezier[80](0.1,0.5)(0.6,0.5)(0.6,0)
%\qbezier[80](0.1,-0.5)(0.6,-0.5)(0.6,0)
%\put(0.6,0){\line(1,0){0.5}}
%\put(0.1,-0.5){\line(0,1){1}}
%\put(-0.4,0.25){\line(1,0){0.5}}
%\put(-0.4,0){\line(1,0){0.5}}
%\put(-0.4,-0.25){\line(1,0){0.5}}
%\end{picture}
\vspace{4pt}
$$
\caption{Elements of $\Lambda$}\label{f:tx3x5}
\end{figure}

The space of primitive elements $\P$ of $\Ab$ becomes a
$\Lambda$-module by
inserting an element of $\Lambda$ into a freely chosen trivalent vertex
of a diagram of a primitive element.
Multiplication by $t$ increases the degree by $1$
and multiplication by $x_3$ 
%($i=3,5$) 
increases the degree by~$3$.
An example is shown in Figure~\ref{f:lmod}.

\begin{figure}[!ht]
$$
x_3\omega_4=
x_3\
\Picture{
\FullCircle
\thinlines
\put(-0.866,0.5){\line(1,0){0.966}}
\put(-0.866,-0.5){\line(1,0){0.966}}
\qbezier[80](0.1,0.5)(0.6,0.5)(0.6,0)
\qbezier[80](0.1,-0.5)(0.6,-0.5)(0.6,0)
\put(0.6,0){\line(1,0){0.4}}
\put(-1,0){\line(1,0){0.5}}
\put(-0.5,-0.5){\line(0,1){1}}
}\quad =\quad
\Picture{
\FullCircle
\thinlines
\put(-0.866,0.5){\line(1,0){0.966}}
\put(-0.866,-0.5){\line(1,0){0.966}}
\qbezier[80](0.1,0.5)(0.6,0.5)(0.6,0)
\qbezier[80](0.1,-0.5)(0.6,-0.5)(0.6,0)
\put(0.6,0){\line(1,0){0.4}}
\put(-1,0){\line(1,0){0.5}}
\put(-0.5,-0.5){\line(0,1){1}}
\put(-0.2,-0.5){\line(0,1){1}}
\put(0.1,-0.5){\line(0,1){1}}
%\put(-0.2,0.25){\line(1,0){0.3}}
\put(-0.2,0){\line(1,0){0.3}}
%\put(-0.2,-0.25){\line(1,0){0.3}}
}\vspace{8pt}
$$
\caption{How $\P$ becomes a $\Lambda$--module}\label{f:lmod}
\end{figure}

If $D$ and $D'$ are classes of trivalent diagrams with a distinguished
oriented circle modulo (STU)-relations (see Figure~\ref{f:rels}),
then their connected sum $D\# D'$ along these circles is well defined.
We state in the following lemma
how the weight systems~$\Wgl$ and~$\Wso$ behave under the
operations described above:
Part~(1) of the lemma is easy to prove; for Part~(2), 
see Theorem~6.4 and Theorem~6.7 of \cite{Vog}.

\begin{lemma}\label{diag_op}
(1) Let $D$ and $D'$ be chord diagrams each one having a distinguished oriented circle.
Then the connected sum of~$D$ and~$D'$ satisfies

$$
\Wgl(D\# D')=\Wgl(D) \Wgl(D')/c \quad \mbox{and} \quad
\Wso(D\# D')=\Wso(D) \Wso(D')/c.
$$

(2) For a primitive element $p\in \P$ we have:

$$
\begin{array}{l}
\Wgl(t p) = c \Wgl(p),\\
\Wso(t p) = \wc \Wso(p),\\
\Wgl(x_3 p) = (c^3+12c) \Wgl(p),\\
\Wso(x_3 p) = (\wc^3-3\wc^2+30\wc-24) \Wso(p),%\\
%\Wgl(x_5 p) = (c^5+32c^3+48c) \Wgl(p),\\
%\Wso(x_5 p) = (\wc^5-5\wc^4+80\wc^3-184\wc^2+408\wc-288)
%\Wso(p),
\end{array}
$$

where $\wc=c-2$.
\end{lemma}

We have the following formulas concerning spaces of
weight systems restricted to primitive elements. 

\begin{prop}\label{dimtab2}
For the restrictions of the weight systems to primitive elements 
of degree $n\geq 1$ we have

\begin{eqnarray*}
(1) & & \dim \Hn1_{\vert\P_n} = \dim \Hn1=%\left[\frac{n}{2}\right]
[n/2],\label{dimhn} \\
%\mbox{for all $n\geq 1$,}\\
(2) & & \dim \Kn1_{\vert\P_n} =\max(n-2, [n/2])=
\left\{
\begin{array}{ll}
%\left[\frac{n}{2}\right] 
[n/2] & \mbox{if $n\leq 3$,}\\
n-2 & \mbox{if $n\geq 3$,}
\end{array}
\right.\\
(3) & & \dim \left(\Hn1_{\vert\P_n} \cap \Kn1_{\vert\P_n}\right) =\min(2, [n/2]) =
\left\{
\begin{array}{ll}
%\left[\frac{n}{2}\right] 
[n/2] & \mbox{if $n\leq 3$,}\\
2 & \mbox{if $n\geq 4$.}
\end{array}
\right.
\end{eqnarray*}
\end{prop}

The proof of Proposition~\ref{dimtab2} will be given in this section together with a proof of 
Theorem~\ref{t:dimtab}. The proof is divided into several steps.

If $q$ is a polynomial, then we denote the degree of its lowest
degree term by
$\ord(q)$.
Now we start to derive lower bounds for dimensions of
spaces of weight systems.

\medskip

\begin{proof}[Proof of Part~(1) of Proposition~\ref{dimtab2}]
By Lemma~\ref{valofdiags} we have $\ord(\Wgl(\omega_{k}))=1$ for even~$k$.
By Lemma~\ref{diag_op} we have
$\Wgl(t^k p)=c^k \Wgl(p)$ for $p\in\P$.
This implies

$$
\dim\left(\Wgl\left(\vspan\{t^{n-2}\omega_2,t^{n-4}\omega_4,
\ldots,t^{n-2[n/2]}\omega_{2[n/2]}\}\right)\right)=
[n/2]\leq \dim\Hn1_{\vert\P_n}\kern-2pt.
$$

Since this lower bound coincides with the upper bound from
Lemma~\ref{l:upperbo} we have $\dim \Hn1_{\vert\P_n}=[n/2]$.
\end{proof}%$\Box$

\medskip

%Define $\omega_n^l=\omega_n\amalg {S^1}^{\amalg (l-1)}$.
Let $D_{ijk}=(L_i\# C_j)\# T_k$ (in this definition we choose arbitrary
distinguished circles of
$L_i$, $C_j$, $(L_i\# C_j)$ and for further use also for
$D_{ijk}$). Let $d_{ijk}$ be the number of oriented circles in $D_{ijk}$
and define $D_{i,j,k}^\ell=D_{ijk}\amalg {S^1}^{\amalg(\ell-d_{ijk})}$
for $\ell\geq d_{ijk}$.
We will make use of the formulas for
$\Wgl(D_{i,0,k}^\ell\#\omega_m)$ and
$\Wso(D_{i,j,k}^\ell\#\omega_m)$ implied by
Lemmas~\ref{valofdiags}
and~\ref{diag_op} throughout the rest of this section.

\medskip

\begin{proof}[Proof of Part~(1) of Theorem~\ref{t:dimtab}]
For all $n\geq 1$ we have $[n/2]$ primitive elements~$p_i$ 
such that the polynomials
$g_i=\Wgl(p_i\amalg {S^1}^{\amalg (\ell-1)})$ are linearly
independent and~$c^\ell\vert g_i$
(see the proof of Part~(1) of Proposition~\ref{dimtab2}).
Let $n<\ell$. The diagrams

\begin{equation}
D_{n,0,0}^\ell, D_{n-2,0,0}^\ell\#\omega_2,\ldots,
D_{n-2[(n-1)/2],0,0}^\ell\#\omega_{2[(n-1)/2]}\label{gllist1}
\end{equation}

are mapped by $\Wgl$ to the values

$$
c^{\ell-n}(1-c^2)^n,c^{\ell-n+2}f_2(c),\ldots,c^{\ell-1}f_{[(n+1)/2]}(c)
$$

with polynomials $f_i$ satisfying $f_i(0)=-2$
($i=2,\ldots,[(n+1)/2]$). So in this case we have found
$
[n/2]+[(n+1)/2]=n
$
linearly independent values, which is the maximal
possible number (see Lemma~\ref{l:upperbo}).
If~$n\geq \ell$, then we conclude in the same way using
the following list of~$k-n+1+[(n-1)/2]=k-[n/2]$ elements where $k=[(n+\ell-1)/2]$:

\begin{equation}
D_{2k-n,0,0}^\ell\#\omega_{2n-2k},
D_{2k-n-2,0,0}^\ell\#\omega_{2n-2k+2},
\ldots, D_{n-2[(n-1)/2],0,0}^\ell\#\omega_{2[(n-1)/2]}.\label{gllist2} %\qed
\end{equation}
\end{proof}%$\Box$

\smallskip

We will use the
upper bounds for~$\dim \Hnl$ and~$\dim\Fnl$ together with the following
lower bound for $\dim (\Hnl\cap\Fnl)$ to get an upper bound for
$\dim (\Hnl+\Fnl)$. In the case~$\ell=1$
we will argue in a similar way for the restriction
of weight systems to primitive elements.

\begin{lemma}\label{capdim}
For all $n,\ell\geq 1$ we have

\begin{eqnarray*}
\dim\left(\Hnl\cap\Fnl\right) & \geq & \min(\dim\Hnl, 2)\\
=\;\min(n,[(n-1+\ell)/2], 2) & = & \dim(\vspan\{W(r_n^\ell),W(y_n^\ell)\}).
\end{eqnarray*}

%\qquad\mbox{and}\qquad
For all $n\geq 1$ we have 

$$
\dim \left(\Hn1_{\vert\P_n}\cap \Fn1_{\vert\P_n}\right)\geq \min(\dim\Hn1,2)=\min([n/2],2).
$$
\end{lemma}
\begin{proof}%{\bf Proof:}
%Let $\ell\geq 1$ and $n\geq 4$. Then we have
%$\dim\Hnl\geq 2$ by Part~(1) of Proposition~\ref{dimtab}.
Propositions~\ref{p:hkpolys} and~\ref{WglWso} imply
that the weight system $W(r_n^\ell)\in\Hnl\cap\Fnl$ is
equal to~$(-1)^\ell\Wgl(.)(2)_{\vert\Abnl}$ and the weight system
$W(y_n^\ell)\in\Hnl\cap\Fnl$ is
equal to the coefficient of~$c^{\ell+n}$ in $\Wgl_{\vert\Abnl}$.
By the proof of Lemma~\ref{l:upperbo}
we have $\Wgl(D)(0)=\Wgl(D)(1)=0$ and
in the weight system $\Wgl_{\vert\Abnl}$
the coefficients of $c^{\ell+n-1},c^{\ell+n-3},\ldots$ and the coefficients of 
$c^{\ell-n-1}, c^{\ell-n-2},\ldots$
vanish. By Part~(1) of
Theorem~\ref{t:dimtab} 
these are the only linear dependencies between
the coefficients of $c^{\ell+n},c^{\ell+n-1},\ldots$ in 
the polynomial~$\Wgl_{\vert\Abnl}$.
This %together with $\dim\Hnl\geq 2$ 
implies for~$\dim\Hnl =1$ that the coefficient of $c^{\ell+n}$ in~$\Wgl_{\vert\Abnl}$ 
is not the trivial weight system 
and this implies for $\dim\Hnl\geq 2$ 
that $\Wgl(.)(2)_{\vert\Abnl}$ and the coefficient of
$c^{\ell+n}$ in $\Wgl_{\vert\Abnl}$ 
are linearly independent.
By Part~(1) of Proposition~\ref{dimtab2} we can argue in the same way
with $\Wgl_{\vert\Pn}$. This completes the proof.
\end{proof}%$\Box$                           

\medskip

Define the weight system $w=(-2)^n\Wso(\cdot)(4)-2(-2)^\ell\Wso(\cdot)(-2)\in\Fnl$.
For~$n\geq 4$ Lemmas~\ref{valofdiags} and~\ref{diag_op} imply that

\begin{equation}
w\left(\omega_2\#(t^{n-4}\omega_2)\amalg {S^1}^{\amalg\ell-1}\right) 
= 18(-4)^n4^{\ell-1}\not=0.
\end{equation}

Part~(3) of Proposition~\ref{p:hkpolys} together with 
Propositions~\ref{WglWsoIntegral} 
and~\ref{WglWso} implies that

\begin{equation}\label{e:wnot0}
0\not=w\in\bigoplus_{i=1}^{n-1}\Fil\Fnmil.
\end{equation}

For~$\ell=1$ Equation~(\ref{e:wnot0}) implies~$w(\P_n)=0$. Therefore we have

\begin{equation}\label{e:betterbo}
\dim\Fn1_{\vert\P_n}\leq \dim\Fn1-1\leq n-2\quad\mbox{for all $n\geq 4$.}
\end{equation}

Since we have $\dim\Hn1=\dim \Hn1_{\vert\P_n}$
by Part~(1) of Proposition~\ref{dimtab2} we 
know that~$w\not\in \Hn1$ and therefore

\begin{equation}\label{e:betterbo2}
\dim (\Hn1+\Fn1)\geq \dim(\Hn1+\Fn1)_{\vert\P_n}+1\quad\mbox{for all $n\geq 4$.}
\end{equation}

Let $(\Wgl,\Wso):\Abnl\longrightarrow \Q[c]\times\Q[c]$ be defined by

$$
\left(\Wgl,\Wso\right)(D)=\left(\Wgl(D),\Wso(D)\right).
$$

Then by Proposition~\ref{WglWso} we have

\begin{eqnarray}
 & \dim(\Hnl+\Fnl) = \dim\left((\Wgl,\Wso)(\Abnl)\right), &
\label{WgloWso}\\
& \dim\left(\Hn1_{\vert\Pn}+\Fn1_{\vert\Pn}\right)=
\dim\left((\Hn1+\Fn1)_{\vert\Pn}\right) =
\dim\left((\Wgl,\Wso)(\Pn)\right).
\label{WgloWsopr}
\end{eqnarray}

We will use Equations~(\ref{WgloWso}) and~(\ref{WgloWsopr})
to derive lower bounds 
for $\dim(\Hnl+\Fnl)$ and 
for~$\dim(\Hn1+\Fn1)_{\vert\Pn}$.
Now we can complete the proofs of Theorem~\ref{t:dimtab} and Proposition~\ref{dimtab2}.

\begin{proof}[Proof of Parts~(2) and~(3) of Proposition~\ref{dimtab2} and 
Theorem~\ref{t:dimtab} for $\ell=1$]

Let 

$$
\Sigma_7=\vspan\{\omega_7,t\omega_6,t^2\omega_5,t^3\omega_4,t^5\omega_2,x_3\omega_4\}\subset\P_7.
$$

Define for $n>7$:

$$
\Sigma_n=\left\{
\begin{array}{ll}
t\Sigma_{n-1}+\Q \omega_n & \mbox{if $n$ is odd,}\\
t\Sigma_{n-1}+\Q \omega_n + \Q x_3\omega_{n-3} & \mbox{if $n$ is even.}
\end{array}
\right.
$$

By a calculation using Lemmas~\ref{valofdiags} and~\ref{diag_op} we obtain

$$
\dim\left((\Wgl,\Wso)(\Sigma_7)\right)=6.
$$

In view of the proof of Lemma~\ref{l:upperbo} we can define a
polynomial-valued
weight system by $\widetilde{W}_\so(.)=\Wso(.)/(c(c-1))$.
We used Lemma~\ref{valofdiags} and Lemma~\ref{diag_op} to compute the
degree~$1$ coefficients of the values of~$\Wgl$ and~$\widetilde{W}_\so$
on elements of~$\Sigma_n$ stated in Table~\ref{coeftab}.

\setlength{\extrarowheight}{1pt}
\begin{table}[!ht]
$$
\begin{array}{|l|c|c|c|c|}
\hline
     & t\Sigma_{n-1} &
%\begin{array}{c}t \Sigma_{n-1}\\ \omega_2\#(t^{n-4}\omega_2)\end{array}  &
     \begin{array}{c}\omega_n\\ \mbox{($n$ odd)}\end{array} &
     \begin{array}{c}\omega_n\\ \mbox{($n$ even)}\end{array} &
     \begin{array}{c} x_3\omega_{n-3}\\ \mbox{($n$ even)}\end{array}\\
\hline
\mbox{coeff. of $\wc$ in $\widetilde{W}_\so(\cdot)$:}& 0 & R_n(2) & R_n(2) & -24 R_{n-3}(2)\\
\hline
\mbox{coeff. of $c$ in $\Wgl(\cdot)$:}    & 0 & 0   & -2  & 0\\
\hline
\end{array}
$$
\caption{Degree~1 coefficients of $\Wgl$ and $\widetilde{W}_\so$ on
$\Sigma_n$}\label{coeftab}
\end{table}

By Lemma~\ref{valofdiags} we have $R_k(2)\not=0$ if $k\not=3$. Then,
by Table~\ref{coeftab} and induction, we see that
$\dim (\Wgl,\Wso)(\Sigma_n)=[n/2]+n-4$ for~$n\geq 7$.
By Equation~(\ref{WgloWsopr}),
Lemmas~\ref{l:upperbo} and~\ref{capdim}, and Equation~(\ref{e:betterbo}) we obtain

\begin{eqnarray*}
[n/2]+n-4+2\leq\dim \left(\Hn1_{\vert\P_n}+\Fn1_{\vert\P_n}\right) & + &
\dim \left(\Hn1_{\vert\P_n}\cap\Fn1_{\vert\P_n}\right)=\nonumber\\
=\dim \Hn1_{\vert\P_n}+\dim \Fn1_{\vert\P_n} & \leq & [n/2]+n-2.\label{upbo2}
\end{eqnarray*}

Thus equality must hold.
This implies Parts~(2) and~(3) of Proposition~\ref{dimtab2} for~$n\geq 7$.
By Equation~(\ref{e:betterbo2}) we get~$\dim (\Hn1+\Fn1)\geq [n/2]+n-3$.
Now we see by Lemmas~\ref{l:upperbo} and~\ref{capdim} that

\begin{eqnarray*}
[n/2]+n-3+2 & \leq & \dim \left(\Hn1+\Fn1\right)  + 
\dim \left(\Hn1\cap\Fn1\right)=\nonumber\\
& = & \dim \Hn1+\dim \Fn1 \leq [n/2]+n-1\label{upbo}
\end{eqnarray*}

which implies Part~(2) and
Part~(3) of Theorem~\ref{t:dimtab} for~$n\geq 7$ and~$\ell=1$.

Let~$\psi$ be the element of degree~$6$ shown in
Figure~\ref{f:psis}.

\begin{figure}[!ht]
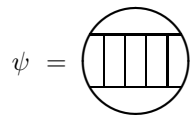

$$
%\psi_1\ = \
%\Picture{
%\FullCircle
%\thinlines
%\put(-0.866,-0.5){\line(1,0){1.732}}
%\put(-0.866,0.5){\line(1,0){1.732}}
%\multiput(-0.625,0.5)(0.25,0){6}{\line(0,-1){1}}
%}\quad , \quad
%\psi_2\ = \
%\Picture{
%\FullCircle
%\thinlines
%\put(-0.866,-0.5){\line(1,0){1.732}}
%\put(-0.866,0.5){\line(1,0){1.732}}
%\put(0.966,0.259){\line(-1,0){0.366}}
%\put(0.966,-0.259){\line(-1,0){0.366}}
%\multiput(-0.6,0.5)(0.4,0){4}{\line(0,-1){1}}
%}\quad , \quad
\psi\ = \
\Picture{
\FullCircle
\thinlines
\put(-0.866,-0.5){\line(1,0){1.732}}
\put(-0.866,0.5){\line(1,0){1.732}}
\multiput(-0.6,0.5)(0.4,0){4}{\line(0,-1){1}}
%}\quad , \quad
%\psi\ = \
%\Picture{
%\FullCircle
%\thinlines
%\put(-0.866,-0.5){\line(1,0){1.732}}
%\put(-0.866,0.5){\line(1,0){1.732}}
%\put(0.966,0.259){\line(-1,0){0.366}}
%\put(0.966,-0.259){\line(-1,0){0.366}}
%\multiput(-0.6,0.5)(0.6,0){3}{\line(0,-1){1}}
}\vspace{8pt}
$$
\caption{A primitive element in degree 6}\label{f:psis}
\end{figure}

A calculation done by computer yields

\begin{eqnarray*}
\Wgl(\psi) & = & c^7+13c^5-14c^3,\\
\widetilde{W}_\so(\psi) & = & \wc^5-3\wc^4+34\wc^3-36\wc^2+16\wc.
%\Wgl(\psi_2) & = & c^9+6c^7+7c^5-14c^3,\\
%\widetilde{W}_\so(\psi_2) & = & \wc^7-6\wc^6+17\wc^5
%                    -25\wc^4+82\wc^3-92\wc^2+48\wc.
\end{eqnarray*}

Let $\Sigma_4=\vspan\{\omega_4,t^2\omega_2\}$,
$\Sigma_5=t\Sigma_4+\Q\omega_5$, and
$\Sigma_6=t\Sigma_5+\Q\omega_6+\Q\psi$.
We obtain again $\dim (\Wgl,\Wso)(\Sigma_n)=[n/2]+n-4$ which implies 
Parts~(2) and~(3) of Proposition~\ref{dimtab2} and 
Theorem~\ref{t:dimtab} for~$\ell=1$ and~$n\geq 4$ 
by the same argument as before. In degrees~$n=1,2,3$ we have 
$\dim\;\P_n=\dim\;\Ab_n=\dim\Hn1=\dim\Fn1=[n/2]$. This completes the proof. 
\end{proof}%$\Box$

\medskip

\begin{proof}[Proof of Parts~(2) and~(3) of Theorem~\ref{t:dimtab} for 
$\ell>1$]
Let $n\geq 4$ and $\ell>1$. 
By the previous proof we have $n+[n/2]-3$ elements
$a_i\in\Ab_n$ such that the values

$$
\left(\Wgl(D_i),\Wso(D_i)\right)\in\Q[c]\times\Q[c]
$$

of $D_i=a_i\amalg {S^1}^{\amalg\ell-1}$ are linearly independent.
%(see Proposition~\ref{dimtab2}).
Consider the following lists of elements:

\begin{eqnarray}
&  & \mbox{If $n\leq \ell$, then we take the $n$ elements}\nonumber\\
& & D_{0,n,0}^\ell,D_{1,n-1,0}^\ell,\ldots,D_{n-3,3,0}^\ell,
D_{0,0,n}^\ell, E_n^\ell:=D_{0,0,n}^\ell-D_{0,0,n-2}^\ell\#\omega_2.
\label{solist1}\\
&  & \mbox{If $n \geq \ell+1$, then we take the $\ell$
elements}\nonumber\\
& & D_{0,\ell-1,n-\ell+1}^\ell,D_{1,\ell-2,n-\ell+1}^\ell,\ldots,
D_{\ell-3,2,n-\ell+1}^\ell,
D_{0,0,n}^\ell, E_n^\ell.\label{solist2}
%& (3) & \mbox{If $n\geq l+2$, then we take the $l$ elements}\\
%& & D_{0,l,0}^l\# \omega_{n-l},D_{1,l-1,0}^l\# \omega_{n-l},\ldots,
%D_{l-2,2,0}^l\#\omega_{n-l},
%D_{0,0,n}^l.
\end{eqnarray}

Let 
${\cal M}_{n,\ell}$ be the list of elements 
$D_i$ together with the elements from
Equation~(\ref{gllist1}) (resp. (\ref{gllist2}))
and Equation~(\ref{solist1}) (resp. (\ref{solist2})).
We have

\begin{equation}\label{cardMnl}
\card\left({\cal M}_{n,\ell}\right) =\left\{
\begin{array}{ll}
3n-3 & \mbox{if $n<\ell$,}\\
n+\ell-3+[(n+\ell-1)/2] & \mbox{if $n\geq \ell$.}
\end{array}
\right.
\end{equation}

The values of $\Wgl$ and $\Wso$ on elements of ${\cal M}_{n,\ell}$ have
the properties stated in Table~\ref{t:Wglsoprop}.

\setlength{\extrarowheight}{2pt}
\begin{table}[!ht]
$$
\begin{array}{|l|l|}
\hline
\ord(\Wgl(D_i))\geq \ell &
\ord(\Wso(D_i))\geq \ell, \Wso(D_i)(2)=0\\
\hline
\ord(\Wgl(E_n^\ell))\geq \ell &
\ord(\Wso(E_n^\ell))\geq \ell, \Wso(E_n^\ell)(2)\not=0\\
\hline
\ord(\Wgl(D_{i,0,0}^\ell\#\omega_{n-i})) & \\
\qquad =\ell-i & \\
\mbox{($i>0$, $n-i$ even)} &
\raisebox{1.6ex}[0cm][0cm]{$\ord(\Wso(D_{i,0,0}^\ell\#\omega_{n-i}))
\geq \ell$}\\
\hline
& \ord(\Wso(D_{n-i,i,0}^\ell))=\ell+1-i\ (i\geq 3)\\
& \ord(\Wso(D_{i,\ell-1-i,n-\ell+1}^\ell))=i+1\ (i\leq \ell-3)\\
& \ord(\Wso(D_{0,0,n}^\ell))=\ell-1\\%, \Wso(D_{0,0,n}^\ell)(2)\not=0\\
\hline
\end{array}
$$
\caption{Properties of $\Wgl(e)$ and $\Wso(e)$ for
$e\in{\cal M}_{n,\ell}$}\label{t:Wglsoprop}
\end{table}

The statements from this table are easily verified. For example, we have

$$
\Wso(E_n^\ell)=c^{\ell-1}(c-1) h(c)
$$

with $h(c)=Q_n(c)-2(c-1)(c-2)Q_{n-2}(c)$.
We have $h(0)=Q_n(0)-4Q_{n-2}(0)=0$ which implies

$$
\ord\left(\Wso(E_n^\ell)\right)\geq \ell,
$$

and $h(2)=Q_n(2)=(-2)^n$ which implies $\Wso(E_n^\ell)(2)\not= 0$.
Now let

$$
f=\sum_{e\in{\cal M}_{n,\ell}} \lambda(e) \left(\Wgl(e),\Wso(e)\right)=
(f_1,f_2)\in \Q[c]\times\Q[c]
$$

be a linear combination with $\lambda(e)\in\Q$.
We want to show that~$f=0$ implies that all scalars
$\lambda(e)$ are~$0$.
For our arguments we will use the entries of Table~\ref{t:Wglsoprop}
beginning at its bottom.
The coefficients $\lambda(D_{n-i,i,0}^\ell)$
(resp. $\lambda(D_{i,\ell-1-i,n-\ell+1}^\ell)$) and
$\lambda(D_{0,0,n}^\ell)$ are~$0$
because they are multiples of

$$
\frac{d^k f_2}{dc^k}(0),\ldots,\frac{d^{\ell-1}f_2}{dc^{\ell-1}}(0)
$$

with $k=\max\{1,\ell-n+1\}$. The coefficients
$\lambda(D_{i,0,0}^\ell\#\omega_{n-i})$ must be $0$ by a similiar argument
for $f_1$. We get $\lambda(E_n^\ell)=0$ because $\Wso(D_i)(2)=0$ and
$\Wso(E_n^\ell)(2)\not=0$. The remaining coefficients $\lambda(D_i)$ are $0$
because the values $(\Wgl(D_i),\Wso(D_i))$ are linearly independent.
This implies $\dim(\Hnl + \Fnl)\geq \card({\cal M}_{n,\ell})$.
By Lemma~\ref{l:upperbo} and Lemma~\ref{capdim} we have

\begin{eqnarray}\label{uneq2}
\card({\cal M}_{n,l})+2 & \leq & \dim\left(\Hnl+\Fnl\right)
+ \dim\left(\Hnl\cap\Fnl\right)= \nonumber\\
=\dim\Hnl+\dim\Fnl & \leq & \left\{
\begin{array}{ll}
3n-1 & \mbox{if $n<\ell$,}\\
n+\ell-1+[(n+\ell-1)/2] & \mbox{if $n\geq \ell$.}
\end{array}
\right.
\end{eqnarray}

Comparing with Equation~(\ref{cardMnl}) shows that equality must hold
in Equation~(\ref{uneq2}).
This completes the proof of Parts~(2) and~(3) of
the theorem for all $n\geq 4$.

In degrees $n=1,2,3$ we used the diagrams shown in Table~(\ref{lowdeg})
(possibly together with some additional circles~$S^1$)
to determine $\dim\Fnl$.

\setlength{\extrarowheight}{0pt}
\begin{table}[!ht]
$$
\begin{array}{|l|l|}
\hline
n=1 & L_1\\
\hline
n=2 & \omega_2, C_2, L_2\\
\hline
n=3 & \omega_3, 
\Omega_3:=\begin{picture}(3.2,0.9)(-1.6,-0.2)
\thicklines
\put(-1,0){\circle{1}}
\put(1,0){\circle{1}}
\thinlines
\put(-0.568,0.27){\line(1,0){1.136}}
\put(-0.5, 0){\line(-1,0){1}}
\put(-0.568,-0.27){\line(1,0){1.136}}
\end{picture}, \omega_2\# L_1, T_3,C_3\\[4pt] 
\hline
\end{array}
$$
\caption{Diagrams used in low degrees}\label{lowdeg}
\end{table}

In the calculation we used the explicit formulas for the values of~$\Wso$
from the proof of Lemma~\ref{valofdiags} together with
$\Wso(\Omega_3)=2c(c-1)(2-c)$.
The number of linearly independent values coincides in all of these cases with the
upper bound for~$\dim \Fnl$ from Lemma~\ref{l:upperbo} or with~$\dim \Abnl$.
For~$\ell\geq 4$ and~$a\in\Ab_{3,3}$ we have

$$
\ord \left(\Wgl\left(L_3\amalg {S^1}^{\amalg \ell-4}\right)\right)=\ell-3\quad\mbox{ and }\quad 
\ord \left(\Wgl\left(a\amalg {S^1}^{\amalg \ell-3}\right)\right)\geq \ell-2.
$$

Together with Lemmas~\ref{l:upperbo} and~\ref{capdim} this implies

$$
3+5-2\geq \dim ({\cal H}'_{3,\ell} + {\cal F}'_{3,\ell})\geq \dim {\cal F}'_{3,3}+1=6
$$

and therefore $\dim ({\cal H}'_{3,\ell} + {\cal F}'_{3,\ell})=6$ and
$\dim ({\cal H}'_{3,\ell} \cap {\cal F}'_{3,\ell})=2$ for~$\ell\geq 4$.
In the cases
$n=1,2$, and in the case $n=3$ and~$\ell<4$, we have~$\dim \Hnl\leq 2$ and obtain
$\dim (\Hnl \cap \Fnl)=\dim \Hnl$
by applying Lemma~\ref{capdim}.
This completes the proof.
\end{proof}%$\Box$

Corollary~\ref{c:hkjy} can now be proven easily.

\begin{proof}[Proof of Corollary~\ref{c:hkjy}]
Proposition~\ref{WglWsoIntegral}, Lemma~\ref{capdim}, and Part~$(3)$ of 
Theorem~\ref{t:dimtab} imply

\begin{eqnarray*}
& & \dim(\vspan\{r_n^\ell, y_n^\ell\})=\dim(\vspan\{W(r_n^\ell), 
W(y_n^\ell)\})\\
& & =\min(\dim\Hnl,2)=\min(\dim\Hvnl,2)=\dim(\Hvnl\cap\Fvnl).
\end{eqnarray*}

By Proposition~\ref{p:hkpolys} we have $r_n^\ell, y_n^\ell\in\Hvnl\cap\Fvnl$. 
This
implies the statement $\vspan\{r_n^\ell, y_n^\ell\}=\Hvnl\cap\Fvnl$ of the
corollary.
\end{proof}

Using Theorem~\ref{t:dimtab} and Proposition~\ref{dimtab2} we can
also prove Theorem~\ref{t:alggen}.

\medskip

\begin{proof}[Proof of Theorem~\ref{t:alggen}]
By Proposition~\ref{dimtab2} we have

$$
\dim(\Hn1+\Fn1)_{\vert\P_n}=b_n:=\left\{\begin{array}{ll}
[n/2] & n\leq 3\\
n+[n/2]-4 & n\geq 4.\\
\end{array}\right.
$$

This implies that in the graded 
algebra~$A$ generated by $\bigoplus_{n=0}^\infty(\Hn1+\Fn1)$ we
find a subalgebra~$B\subseteq A$ which is a polynomial algebra with~$b_n$ 
generators in degree~$n$.
For~$n\geq 4$ we find by Equation~(\ref{e:wnot0}) a nontrivial element~$w\in\Fn1$ 
lying in the algebra generated by $\bigoplus_{n=1}^{n-1}\Fn1$. This shows that~$A$ is
generated by~$a_n$ elements in degree~$n$ with~$a_n:=\dim(\Hn1+\Fn1)-1$ for~$n\geq 4$ 
and~$a_n:=\dim(\Hn1+\Fn1)$ for~$n\leq 3$. By Theorem~\ref{t:dimtab} we have~$a_n=b_n$.
Now~$B\subseteq A$ implies~$A=B$. By Proposition~\ref{WglWsoIntegral}
the isomorphism~$Z^*$ maps~$A$ to the algebra generated 
by~$\bigoplus_{n=0}^\infty(\Hvn1+\Fvn1)$. This completes the proof.
\end{proof}

%\input hkb.tex

%\bigskip

%Jens Lieberum

%Mathematisches Institut

%Sidlerstr. 5

%3012 Bern

%Switzerland

%Email: lieberum@math-stat.unibe.ch

}

%--------------------Here the manuscript ends-------------------------------
\Addresses

\begin{thebibliography}{HOM}\addcontentsline{toc}{section}{\numberline{}Bibliography}
\itemsep=0pt
\bibitem[BN1]{BN1} 
D.\ Bar-Natan, {\em On the Vassiliev knot invariants}, Topology 34
(1995), 423--472.

\bibitem[BNG]{BNG} D.\ Bar-Natan and S.\ Garoufalidis, {\em On the
Melvin-Morton-Rozansky conjecture},
Invent.\ Math.\ 125 (1996), 103--133.

\bibitem[Bra]{Bra} R.\ Brauer, {\em On algebras which are connected with the
semisimple continous groups}, Annals of Math. 38 (1937), 857--872.

%\bibitem[ChD]{ChD}
%S. Chmutov, S. Duzhin, {\em A lower bound for the number of Vassiliev knot
%invariants}, preprint, 1996.

\bibitem[CoG]{CoG} R.\ Correale, E.\ Guadagnini, {\em Large-N Chern-Simons
field theory}, Phys. Letters B 337 (1994), 80--85.

\bibitem[Jon]{Jon} V.\ F.\ R.\ Jones, {\em A polynomial invariant for
knots via von Neumann algebras}, Bull. of AMS 12 (1985), 103 --111.

\bibitem[HOM]{HOM} P.\ Freyd, J.\ Hoste, W.\ B.\ R.\ Lickorish,
K.\ Millet, A.\ Ocneanu and D.\ Yetter,
{\em A new polynomial invariant of knots and links},
Bull. of AMS 12 (1985), 239--246.

\bibitem[KMT]{KMT} T.\ Kanenobu, Y.\ Miyazawa and A.\ Tani, {\em Vassiliev link invariants
of order three}, J.\ of Knot Theory and its Ramif. 7, No.\ 4 (1998), 433--462.

\bibitem[Kas]{Kas} C.\ Kassel, {\em Quantum groups},
GTM 155, Springer-Verlag, New York 1995.

\bibitem[KaT]{KaT} C.\ Kassel and V.\ Turaev, {\em Chord diagram 
invariants of tangles and
graphs}, Duke Math.\ J.\ 92 (1998), 497--552.

\bibitem[Ka1]{Ka1} L.\ H.\ Kauffman, {\em State models and the Jones polynomial},
Topology~26 (1987), 395--407.

\bibitem[Ka2]{Ka2} L.\ H.\ Kauffman, {\em An invariant of regular isotopy}, Trans.\ Am.\
Math.\ Soc.\ 318, No.\ 2 (1990), 417--471.


\bibitem[Ka3]{Kau} L.\ H.\ Kauffman, {\em Knots and Physics (second edition)}, 
World Scientific, Singapore 1993.

\bibitem[Kn1]{Kn1} J.\ A.\ Kneissler, {\em Woven braids and their closures},
J.\ of Knot Theory and its
Ramif.\ 8, No.~2 (1999), 201--214.

\bibitem[Kn2]{Kn2} J.\ A.\ Kneissler, {\em The number of primitive
Vassiliev invariants up to degree twelve}, q-alg/9706022 and
University of Bonn preprint, June 1997.

\bibitem[Lam]{Lam} H.\ Lamaugarny, {\em Sp\'{e}cialisations communes entre
le polyn\^{o}me de Kauffman et le polyn\^{o}me de Jones--Conway},
C.\ R.\ Acad.\ Sci.\ Paris, t.\ 313, S\'{e}rie I (1991), 289--292.

\bibitem[LM1]{LM1} T.\ Q.\ T.\ Le and J.\ Murakami, {\em Kontsevich's 
integral for the Homfly polynomial
and relations between values of multiple zeta functions}, Topology 
and its Appl.~62 (1995), 193--206.

\bibitem[LM2]{LM2} T.\ Q.\ T.\ Le and J.\ Murakami, {\em Kontsevich integral for the Kauffman
polynomial}, Nagoya Math.\ J., 142 (1996), 39--65.

\bibitem[LM3]{LM3} T.\ Q.\ T.\ Le and J.\ Murakami, {\em The universal Vassiliev-Kontsevich invariant
for framed oriented links}, Comp. Math.\ 102 (1996), 41--64.

\bibitem[Lie]{Lie} J.\ Lieberum, {\em Chromatic weight systems and the
corresponding knot invariants}, 
Math.\ Ann., to appear.

\bibitem[Lik]{Lik} W.\ B.\ R.\ Lickorish, 
{\em Polynomials for links}, London Math.\ Soc.\ Bulletin 20 (1988), 558--588.

\bibitem[Men]{Men} G.\ Meng, {\em Bracket models for weight systems and
the universal Vassiliev invariants}, Topo\-logy and its Appl. 76 (1997), 47-60.

\bibitem[Sul]{Sul} P.\ Sulpice, {\em Invariants de n\oe uds et filtration
de Vassiliev}, Th\`ese de Doctorat de l'Universit\'e Paris~VII, 1997.

\bibitem[Vog]{Vog} P. Vogel,
{\em Algebraic structures on modules of diagrams},
Invent.\ Math., to appear.
\end{thebibliography}
\end{document}